\documentclass[mathpazo]{cicp}

\usepackage[dvips]{color}
\newcommand{\halb}{\frac{1}{2}}

\newcommand{\A}{{\mathbf{A}}}
\newcommand{\F}{\mathcal{\bf{F}}} 
\newcommand{\IF}{\textnormal{ if }}

\newcommand{\PNM}{P_NP_M}
\renewcommand{\d}{\partial}
\newcommand{\Q}{\mathbf{Q}}
\renewcommand{\F}{\mathbf{F}}
\renewcommand{\S}{\mathbf{S}}
\newcommand{\w}{\mathbf{w}}
\newcommand{\q}{\mathbf{q}}
\newcommand{\be}{\begin{equation}}
\newcommand{\ee}{\end{equation}}
\newcommand{\bdm}{\begin{displaymath}}
\newcommand{\edm}{\end{displaymath}}
\newcommand{\bea}{\begin{eqnarray}}
\newcommand{\eea}{\end{eqnarray}}
\begin{document}
%%%%% title : short title may not be used but TITLE is required.
% \title{TITLE}
% \title[short title]{TITLE}
\title{On Arbitrary-Lagrangian-Eulerian One--Step WENO Schemes for Stiff Hyperbolic Balance Laws}

 %multiple authors:
 %Note the use of \affil and \affilnum to link names and addresses.
 %The author for correspondence is marked by \corrauth.
 %use \emails to provide email addresses of authors
 %e.g. below example has 3 authors, first author is also the corresponding
 %     author, author 1 and 3 having the same address.
\author[Dumbser et.~al.]{Michael Dumbser\affil{1}\comma\corrauth,
                         Ariunaa Uuriintsetseg\affil{1}, 
                         Olindo Zanotti\affil{1}}  
\address{\affilnum{1}\ Laboratory of Applied Mathematics,
          					   University of Trento,
          						 I-38123 Trento, Italy
         %\affilnum{2}\ Department of Mathematics,
         % 							Hong Kong Baptist University, Hong Kong SAR}
         }
\emails{{\tt michael.dumbser@ing.unitn.it} (M.~Dumbser) }

%%%%% Begin Abstract %%%%%%%%%%%
\begin{abstract}
In this article we present a new family of high order accurate Arbitrary Lagrangian-Eulerian one--step WENO finite volume schemes for 
the solution of stiff hyperbolic balance laws. High order accuracy in space is obtained with a standard WENO reconstruction algorithm 
and high order in time is obtained using the local space-time discontinuous Galerkin method recently proposed in \cite{DumbserEnauxToro}. 
In the Lagrangian framework considered here, the local space--time DG predictor is based on a weak formulation of the governing PDE on 
a moving space--time element. For the space--time basis and test functions we use Lagrange interpolation polynomials defined by 
tensor--product Gauss--Legendre quadrature points. The moving space--time elements are mapped to a reference element using an 
isoparametric approach, i.e. the space--time mapping is defined by the same basis functions as the weak solution of the PDE. 
We show some computational examples in one space--dimension for non--stiff and for stiff balance laws, in particular for the 
Euler equations of compressible gas dynamics, for the resistive relativistic MHD equations, and for the relativistic radiation 
hydrodynamics equations. Numerical convergence results are presented for the stiff case up to sixth order of accuracy in space and 
time and for the non--stiff case up to eighth order of accuracy in space and time.  
\end{abstract}
%%%%% end %%%%%%%%%%%

%%%%% AMS/PACs/Keywords %%%%%%%%%%%
%\pac{}
%\ams{52B10, 65D18, 68U05, 68U07}

%%%% maketitle %%%%%
\maketitle

%%%% Start %%%%%%
\section{Introduction}
\label{sec.intro}

We present a new class of high order one-step Arbitrary Lagrangian-Eulerian (ALE) finite volume schemes for stiff hyperbolic balance 
laws. While the mesh is fixed in an Eulerian description, in Lagrangian type schemes the computational mesh moves with the local fluid 
velocity. That means that material interfaces are moving together with the mesh and thus one can precisely identify their location. 
In the recent past, a lot of work has been carried out to develop Lagrangian methods. Some algorithms are developed 
starting directly from the the conservative quantities such as mass, momentum and total energy \cite{Maire2007, Smith1999} 
while another class starts from the  nonconservative form of the governing equations \cite{Neumann1950,Benson1992,Caramana1998}. 
In any discrete scheme one has to decide where to place the degrees of freedom of each physical variable. The existing Lagrangian 
schemes in literature can be generally separated into two main classes: 1) staggered mesh methods, where the velocity is defined 
at the cell interfaces while the other physical variables are located at the cell center and 2) cell-centered methods, where all 
variables are defined at the cell center. 

In \cite{munz94} Munz presented several different Godunov-type finite volume schemes for Lagrangian gas dynamics and, in particular,  
he was the first to introduce a Roe linearization for Lagrangian gas dynamics. It was found that the Lagrangian 
Roe linearization actually does \textit{not} coincide with the Eulerian one \cite{munz94}. The resulting maximum signal speeds of this Roe 
linearization have subsequently been used to construct robust HLL-type Riemann solvers in Lagrangian coordinates. 
Carr\'e et al. \cite{Carre2009} describe a cell-centered Godunov scheme for Lagrangian gas dynamics on general multi-dimensional 
unstructured meshes. Their finite volume solver is node based and compatible with the mesh displacement. In \cite{Despres2005}, Despr\'es 
and Mazeran propose a way of writing the equations of gas dynamics in Lagrangian coordinates in two dimensions as a weakly hyperbolic 
system of conservation laws. The system contains both, the physical and the geometrical part. Based on the symmetrization of the formulation 
of the physical part, the authors design a finite volume scheme for the discretization of Lagrangian gas dynamics on moving meshes. 
In  \cite{Jua2011}, Jua and Zhang present a high-order Lagrangian Runge-Kutta DG scheme for the discretization of two-dimensional 
compressible gas dynamics. The scheme uses a fully Lagrangian form of the gas dynamics equations and employs a new HWENO-type 
reconstruction algorithm as limiter to control the spurious oscillations, maintaining the compactness of RKDG methods. 
The time marching for the semi discrete schemes for the physical and geometrical variables is implemented by a classical TVD Runge-Kutta 
method. The scheme has been shown to achieve second order of accuracy, both in space and time.  Another Lagrangian discontinuous Galerkin 
finite element method has been recently proposed in \cite{Dimitri2011}. The method preserves discrete conservation in the presence of 
arbitrary mesh motion and thus obeys the Geometric Conservation Law (GCL).

In a series of articles  \cite{Maire2007,Maire2010,Maire2011} Maire et al. develop a general formalism to derive first and second order  cell-centered Lagrangian schemes in multiple space dimensions and also on general polygonal grids. In \cite{Maire2011} the time derivatives 
of the fluxes are obtained through the use of a node-centered solver which can be viewed as a multi-dimensional extension of the 
Generalized Riemann problem methodology introduced by Ben-Artzi and Falcovitz \cite{Artzi}, Le Floch et al.  \cite{Raviart.GRP.1,Raviart.GRP.2} 
and Titarev and Toro \cite{toro3,toro4,titarevtoro}.  In their recent papers \cite{chengshu1,chengshu2} Cheng and Shu developed a class of 
cell centered Lagrangian finite volume schemes for solving the Euler equations which are based on high order essentially non-oscillatory 
(ENO) reconstruction, both with Runge-Kutta and Lax-Wendroff-type timestepping. To our knowledge, the Lagrangian schemes developed in \cite{chengshu1,chengshu2} are the first better than second order non-oscillatory Lagrangian finite volume schemes published so far.  
Further work of Cheng and Shu contains the construction of symmetry-preserving Lagrangian schemes, see \cite{chengshu3,chengshu4}. 

A completely different class of fully Lagrangian methods can be found in meshless particle schemes such as the SPH approach \cite{Monaghan1994,Dambreak3D,SPH3D,SPHWeirFlow,SPHLagrange}, which has become very popular to simulate fluid motion in complex 
deforming domains due to its algorithmic simplicity and high versatility and flexibility. 

Apart from real Lagrangian methods, where the mesh actually moves with the local fluid velocity, and Arbitrary Lagrangian Eulerian (ALE) schemes, see e.g. \cite{Hirt1974,Smith1999,Peery2000}, where the mesh moves with an arbitrary mesh velocity that may or may not coincide with the real  fluid velocity \cite{HuiCoord}, there also exist Semi-Lagrangian schemes. This kind of method is used in general to solve transport equations.  Here, the discrete solution is represented on a fixed Eulerian grid. However, the solution at a mesh point at the new time $t^{n+1}$ is computed from the known solution at time $t^n$ by following back the Lagrangian trajectories of the fluid to the end-point of the trajectory, which in general does not coincide with a grid point. The unknown solution at the end-point of the Lagrangian trajectory is then obtained via interpolation from the known discrete solution at time $t^n$ at surrounding mesh points, see \cite{CIR,Casulli1990,CasulliCheng1992,LentineEtAl2011,QuiShu2011}. 

In this article we introduce a new and better than second order accurate Lagrangian one-step WENO finite volume scheme for the solution of stiff and non-stiff nonlinear systems of 
hyperbolic balance laws. The high order of accuracy in space is obtained using a WENO reconstruction \cite{shu_efficient_weno,balsarashu,DumbserEnauxToro} and the one-step time 
discretization is based on a high order accurate predictor, for which a local space-time discontinuous Galerkin finite element scheme is used \cite{DumbserEnauxToro,DumbserZanotti,HidalgoDumbser}. 

The outline of this article is as follows: in Section
\ref{sec.method} we describe the numerical scheme, while in Section \ref{sec.test} we show numerical results for three
different sets of equations, namely for the compressible Euler equations, for the resistive relativistic MHD equations 
(which provides a natural benchmark of stiff problems) and for the relativistic radiation hydrodynamics equations. 
Results for shock tube problems are shown as well as numerical convergence results for smooth solutions to validate that 
the designed order of accuracy of our schemes is reached. 
Finally, in Section \ref{sec.conclusions} we summarize our results and give some concluding remarks and an outlook concerning possible future extensions of our method. 

\section{Numerical method}
\label{sec.method}

In this article we consider general nonlinear systems of hyperbolic balance laws of the form
\begin{equation}
\label{eqn.pde}
  \frac{\d \Q}{\d t} + \frac{\d \F(\Q)}{\d x} = \S(\Q), \qquad x \in \Omega(t) \subset \mathbb{R}, \quad t \in \mathbb{R}_0^+, \quad \Q \in \Omega_{\Q} \subset \mathbb{R}^\nu,     
\end{equation} 
where $\Q$ is the vector of conserved variables, $\Omega_{\Q}$ is the space of admissible states (state--space), $\Omega(t)$ is the variable
spatial computational domain, $\F(\Q)$ is the flux vector and $\S(\Q)$ is a nonlinear algebraic source term, which can also be stiff. 

\subsection{ALE-Type One--Step Finite Volume Schemes}
\label{sec.fv}
The computational domain $\Omega$ is discretized by a set of \textit{moving} mesh points $x_{i+\halb}=x_{i+\halb}(t)$ that move with a  \textit{general mesh velocity} $V_{i+\halb}=V_{i+\halb}(\Q_{i+\halb},x_{i+\halb},t)$, i.e. 
\begin{equation}
\label{eqn.meshmove} 
\frac{d}{dt} x_{i+\halb} = V_{i+\halb}(\Q_{i+\halb},x_{i+\halb},t). 
\end{equation}  
The spatial control volumes are defined at the current time $t^n$ as $T_i^n = [x_{i-\halb}^n;x_{i+\halb}^n]$, where we have used the notation $x_{i+\halb}^n = x_{i+\halb}(t^n)$. 
By integration over the moving space--time control volume $[x_{i-\halb}(t);x_{i+\halb}(t)] \times [t^n;t^{n+1}]$ and application of Gauss' theorem, the following integral formulation 
for the balance law \eqref{eqn.pde} is obtained: 
\begin{equation}
\label{eqn.fv} 
  \Delta x_i^{n+1} \Q_i^{n+1} = \Delta x_i^{n} \Q_i^{n} - \Delta t \left( \F^V_{i+\halb} - \F^V_{i-\halb} \right) + \Delta x_i^{n} \, \Delta t \S_i,      
\end{equation} 
with the mesh spacing at time $t^n$ $\Delta x_i^{n} = x_{i+\halb}^n  - x_{i-\halb}^n $ and the time step $\Delta t = t^{n+1} - t^n$.  The cell average at time $t^n$ is defined as 
\begin{equation}
  \Q_i^n = \frac{1}{\Delta x_i^{n}} \int \limits_{x_{i-\halb}^n}^{x_{i+\halb}^n} \Q(x,t^n) dx,    
\end{equation}  
the source term is given by   
\begin{equation}
  \S_i = \frac{1}{\Delta x_i^{n} \, \Delta t} \int \limits_{t^n}^{t^{n+1}} \int \limits_{x_{i-\halb}(t)}^{x_{i+\halb}(t)} \S\left(\Q(x,t)\right) \ dx dt,    
\end{equation}  
and the flux at the cell interface on the moving mesh is defined as 
\begin{equation}
  \F^V_{i+\halb} = \frac{1}{\Delta t} \int \limits_{t^n}^{t^{n+1}} \left( \F(x_{i+\halb}(t),t) - V_{i+\halb}(t) \Q(x_{i+\halb}(t),t) \right) \ dt.    
\end{equation}
By choosing $V_{i+\halb}=0$, Eq.~\eqref{eqn.fv} reduces to a classical \textit{Eulerian} finite volume scheme on a fixed mesh, while a   \textit{Lagrangian-type} finite volume scheme is obtained by choosing $V_{i+\halb}$ to be the local fluid velocity. Obviously, any other 
choice of $V_{i+\halb}$ is also possible and within the Arbitrary--Lagrangian--Eulerian (ALE) framework. For convenience of notation, 
we also introduce  
\begin{equation}
\F^V(\Q,V) = \F(\Q) - V \Q, \qquad \textnormal{ and } \qquad \A^V(\Q,V) = \frac{\partial \F^V}{\partial \Q}. 
\end{equation}
While eqn. \eqref{eqn.fv} is an \textit{exact} relation, a numerical scheme is obtained by
using a \textit{numerical flux} $\F^V_h(\q_h^-,\q_h^+)$ , which is a function of \textit{two} arguments, namely the states $\q_h^- = \q_h(x_{i+\halb}^-,t)$ and 
$\q_h^+ = \q_h(x_{i+\halb}^+,t)$ on the left and on the right of the interface, respectively: 
\begin{equation} 
\label{eqn.numflux} 
\F^V_{i+\halb} := \frac{1}{\Delta t} \int \limits_{t^n}^{t^{n+1}} \F^V_h \left(\q_h(x_{i+\halb}^-,t),\q_h(x_{i+\halb}^+,t) \right) dt.  
\end{equation} 
In this article, we use two different numerical fluxes, either a simple Rusanov--type flux \cite{Rusanov:1961a}, or an Osher--type flux 
as introduced in \cite{OsherUniversal}. The Rusanov--type flux reads 
\begin{equation}
\label{eqn.rusanov} 
  \F^V_h(\q_h^-,\q_h^+) = \frac{1}{2} \left( \F^V(\q_h^-,V_{i+\halb}) + \F^V(\q_h^+,V_{i+\halb}) \right) - \frac{1}{2} s_{\max} \left( \q_h^+ - \q_h^- \right),   
\end{equation} 
where $s_{\max} = \max( \max( |\lambda(\A^V(\q_h^-,V_{i+\halb}))|), \max( |\lambda(\A^V(\q_h^+,V_{i+\halb}))|) $ is the maximum signal speed. 
The Osher--type flux according to \cite{OsherUniversal} reads 
\begin{equation}
\label{eqn.osher} 
  \F^V_h(\q_h^-,\q_h^+) = \frac{1}{2} \left( \F^V(\q_h^-,V_{i+\halb}) + \F^V(\q_h^+,V_{i+\halb}) \right) - \frac{1}{2} \left( \int \limits_0^1 \left|\A^V( \mathbf{\Psi}(s), V_{i+\halb})\right| ds \right) \left( \q_h^+ - \q_h^- \right),   
\end{equation} 
where 
\begin{equation}
  \mathbf{\Psi}(s) = \mathbf{\Psi}(\q_h^-,\q_h^+,s) = \q_h^- + s \left(\q_h^+ - \q_h^- \right) 
\end{equation} 
is a straight--line segment path connecting the two states $\q_h^-$ and $\q_h^+$, respectively, and the integral in Eqn. \eqref{eqn.osher} is evaluated \textit{numerically} using appropriate high order Gauss--Legendre quadrature
formulae (see \cite{OsherUniversal} for details).  In \eqref{eqn.osher}, the usual definition for the absolute value of a matrix holds: 
\begin{equation}
 | \A | = \mathbf{R} |\mathbf{\Lambda}| \mathbf{R}^{-1}, 
\end{equation} 
where $\mathbf{R}$ is the matrix of right--eigenvectors, $\mathbf{R}^{-1}$ is its inverse and $|\mathbf{\Lambda}|$ is the diagonal matrix of 
the absolute values of the eigenvalues of $\A$. 

For the mesh velocity, needed in \eqref{eqn.meshmove} and in the fluxes \eqref{eqn.rusanov} and \eqref{eqn.osher} we use the Roe--averaged velocity for Lagrangian gas dynamics 
according to Munz \cite{munz94}:       
\begin{equation}
\label{eqn.meshvelocity}  
V_{i+\halb} = \frac{1}{2} \left( V(\q_h^-) + V(\q_h^+) \right). 
\end{equation}  
Note that in Lagrangian gas dynamics, the Roe average for the velocity is just simply given by the arithmetic average (see \cite{munz94} for details on the derivation) and 
\textit{not} by the common expression  valid in Eulerian
coordinates \cite{roe}. 

Finally, the new position of the mesh point $x_{i+\halb}$ at time $t^{n+1}$ becomes with \eqref{eqn.meshmove} and \eqref{eqn.meshvelocity} 
\begin{equation}
x_{i+\halb}^{n+1} = x_{i+\halb}^n + \int \limits_{t^n}^{t^{n+1}} V_{i+\halb} dt. 
\end{equation} 

Furthermore, in \eqref{eqn.fv} also a discrete form of the source term $\S_i$ must be chosen. Since equation \eqref{eqn.fv} only gives an evolution equation for the cell 
averages $\Q_i^n$ but the interface flux $\F^V_{i+\halb}$ needs values at the element interface, a spatial reconstruction operator is needed that produces suitable interface 
values from the given cell averages. The original first order Godunov finite volume scheme uses the simple reconstruction 
\begin{equation}
\q_h(x_{i+\halb}^-,t) = \Q_i^n, \qquad \textnormal{and} \qquad \q_h(x_{i+\halb}^+,t) = \Q_{i+1}^n.
\end{equation} 
Higher order spatial and temporal accuracy can be obtained by using a more sophisticated reconstruction operator, described in the following section.  

\subsection{Polynomial WENO Reconstruction on Irregular Meshes}
\label{sec.weno}

In this paper we use the \textit{polynomial} WENO reconstruction algorithm proposed in \cite{DumbserKaeser06b,DumbserKaeser07,DumbserEnauxToro} that produces
as output entire reconstruction polynomials and not point values at the cell interfaces, as the original optimal WENO scheme of Jiang and Shu \cite{shu_efficient_weno}. 
Since the details can be found in the above--mentioned references, here we only give a brief summary of the algorithm supposing \textit{componentwise} reconstruction
in \textit{conservative variables}. For more details on reconstruction in \textit{characteristic variables} see \cite{eno,shu_efficient_weno}. 
The reconstruction polynomial of degree $M$ is obtained componentwise by requiring \textit{integral conservation} on a stencil 
\begin{equation}
\label{eqn.stencildef}  
\mathcal{S}_i^s = \bigcup \limits_{j=i-l}^{i+r} T_j^n   
\end{equation}
with spatial extension $l$ and $r$ to the left and right, respectively. For odd order schemes there is only one central stencil ($s=1$), with $l=r=M/2$. For even order schemes,
there are two central stencils with $l=$floor$(M/2)+1$ and $r=$floor$(M/2)$ for the first central stencil ($s=0$) and $l=$floor$(M/2)$, $r=$floor$(M/2)+1$ for the 
second one ($s=1$). For all schemes, the fully left--sided stencil ($s=2$) has $l=M$ and $r=0$ and the fully right--sided stencil has $l=0$ and $r=M$. 
The reconstruction polynomial for each candidate stencil is written in terms of some spatial basis functions $\psi_m(\xi)$ as  
\begin{equation}
\label{eqn.recpolydef} 
 \w^s_h(x,t^n) = \sum \limits_{m=0}^M \psi_m(\xi) \hat \w^s_m := \psi_m(\xi) \hat \w^{n,s}_m,   
\end{equation}
with the mapping to the reference coordinate given by 
\begin{equation}
 x = x(\xi,i) = x_{i-\halb}^n + \Delta x_i^n \xi,  \qquad \textnormal{and} \qquad 
 \xi = \xi(x,i) = \frac{1}{\Delta x_i^n} \left( x-x_{i-\halb}^n \right). 
\end{equation} 
Throughout this paper we use the Einstein summation convention, implying summation over indices appearing twice. We furthermore use the Legendre polynomials rescaled to the 
unit interval $I=[0;1]$ as basis functions $\psi_m(\xi)$. Integral conservation on all elements of the stencil then yields  
\begin{equation}
 \frac{1}{\Delta x^n_j} \int \limits_{x^n_{j-\halb}}^{x^n_{j+\halb}} \psi_m(x) \hat \w^{n,s}_m dx = \Q^n_j, \qquad \forall T_j^n \in \mathcal{S}_i^s.     
\end{equation}
and, with the definition of the primitive functions of $\psi_m(\xi)$
\begin{equation}
\label{eqn.prim} 
\bar \Psi_m(\xi) = \int \limits_{0}^{\xi} \psi_m(\zeta) d\zeta, 
\end{equation}
one obtains the following compact expression for the linear algebraic system that has to be solved for the unknown coefficients $\hat w^{n,s}_m$: 
\begin{equation}
 \left( \bar \Psi_m(\xi(x^n_{j+\halb},i)) - \bar\Psi_m(\xi(x^n_{j-\halb},i)) \right) \hat \w^{n,s}_m = \left( \xi(x^n_{j+\halb},i) - \xi(x^n_{j-\halb},i) \right) \Q^n_j, \qquad \forall T_j^n \in \mathcal{S}_i^s.     
\end{equation}
Adopting the usual definitions of the oscillation indicators $\sigma_s$ \cite{shu_efficient_weno} and the oscillation indicator matrix $\Sigma_{lm}$ \cite{DumbserEnauxToro} 
\begin{equation}
\sigma_s = \Sigma_{lm} \hat \w^{n,s}_l \hat \w^{n,s}_m, \qquad 
\Sigma_{lm} = \sum \limits_{\alpha=1}^M \int \limits_0^1 \frac{\partial^\alpha \psi_l(\xi)}{\partial \xi^\alpha} \cdot \frac{\partial^\alpha \psi_m(\xi)}{\partial \xi^\alpha} d\xi,   
\end{equation} 
the nonlinear weights $\omega_s$ are defined by
\begin{equation}
\tilde{\omega}_s = \frac{\lambda_s}{\left(\sigma_s + \epsilon \right)^r}, \qquad 
\omega_s = \frac{\tilde{\omega}_s}{\sum_q \tilde{\omega}_q},  
\end{equation} 
where we use $\epsilon=10^{-14}$, $r=8$, $\lambda_s=1$ for the one--sided stencils and $\lambda=10^5$ for the central stencils, according to \cite{DumbserEnauxToro,DumbserKaeser06b}. 
The final nonlinear WENO reconstruction polynomial and its coefficients are then given by 
\begin{equation}
\label{eqn.weno} 
 \w_h(x,t^n) = \psi_m(\xi) \hat \w^{n}_m, \qquad \textnormal{ with } \qquad  
 \hat \w^{n}_m = \sum_s \omega_s \hat \w^{n,s}_m.   
\end{equation}   

\subsection{Local Space--Time DG Predictor on Moving Meshes}
\label{sec.stdg}

The reconstruction polynomials $\w_h(x,t^n)$ are then evolved \textit{locally} within each element in order to obtain high order time accuracy. 
Instead of the Cauchy--Kovalewski procedure based on Taylor series and repeated differentiation of the governing PDE used in the original ENO 
method of Harten et al. \cite{eno}, in the ADER schemes
of Titarev and Toro
\cite{toro3,toro4,titarevtoro,titarevtoro2,Toro:2006a}
and in the Lagrangian 
ENO finite volume scheme with Lax--Wendroff time discretization presented by Liu et al. \cite{chengshu2}, we use a \textit{weak formulation} of the 
governing PDE in space--time based on the local space--time discontinuous Galerkin method introduced in  \cite{DumbserEnauxToro,DumbserZanotti,HidalgoDumbser}, 
which is also capable of dealing with \textit{stiff} algebraic source terms. Due to the element--local formulation, the method proposed here is different from 
the global space--time DG method of Van der Vegt and Van der Ven \cite{spacetimedg1,spacetimedg2}. 

\begin{figure}[!htbp]
\begin{center}
\begin{tabular}{lr}
\includegraphics[width=0.45\textwidth]{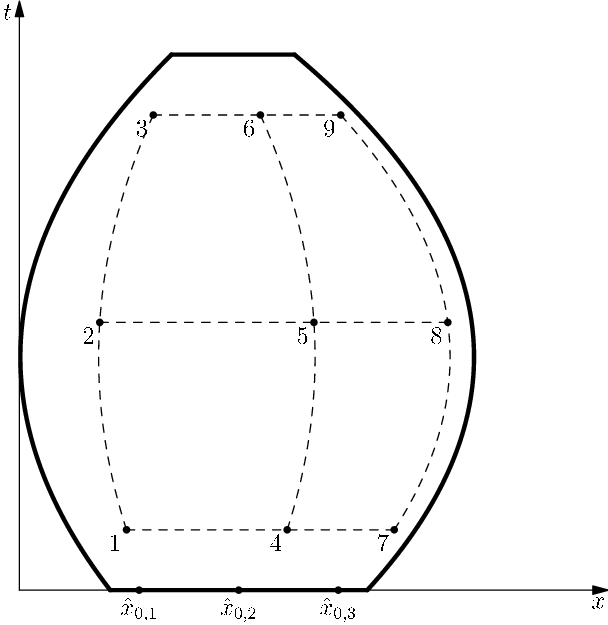} &
\includegraphics[width=0.45\textwidth]{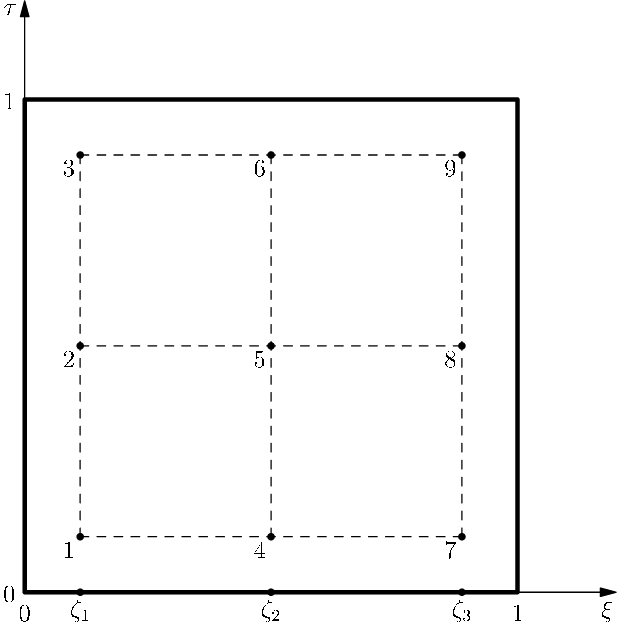}  
\end{tabular}
\caption{Sketch of a third order isoparametric space--time element. Left: physical space--time element. Right: reference space--time element. The interpolation nodes for 
the numerical solution and for the mapping, given by the tensor--product Gauss--Legendre quadrature points, are numbered from 1 to 9. The initial location for the spatial
Gauss--Legendre nodes $\hat x_{0,m}$ is also highlighted. }
\label{fig.isoparam}
\end{center}
\end{figure}

In order to get an element--local weak formulation of the PDE on the \textit{moving} space--time control volume $[x_{i-\halb}(t);x_{i+\halb}(t)]\times[t^n;t^{n+1}]$, the 
governing PDE \eqref{eqn.pde} is first transformed to the reference space-time element $T_E=[0;1]^2$. Therefore, we \textit{map} the physical variables $x$ and $t$ onto 
the reference variables $\xi$ and $\tau$, using an \textit{isoparametric} mapping, i.e. for the mapping of the coordinates we use the same basis functions $\theta_m$ that 
are also used to represent the numerical solution. In this article we use for $\theta_m$ the \textit{Lagrange} interpolation polynomials of degree $M$ that pass through the 
tensor--product \textit{Gauss--Legendre} quadrature points on the reference element $T_E=[0;1]^2$. For details on multidimensional quadrature formulae 
see \cite{stroud}. In the following, the underlying one--dimensional Gauss--Legendre quadrature points and weights on the unit interval $[0;1]$ are denoted by $\zeta_j$ and 
$\alpha_j$, respectively. Using the space-time basis functions $\theta_k$, the mapping of $x$ and $t$ onto $\xi$ and $\tau$ simply reads  
\begin{equation} 
\label{eqn.iso.map} 
 x(\xi,\tau) = \hat x_m \theta_m(\xi,\tau), \qquad  
 t(\xi,\tau) = \hat t_m \theta_m(\xi,\tau).  
\end{equation} 
Here, $\theta_m = \theta_m(\xi,\tau)$ and the coefficients $\hat x_m$ and $\hat t_m$ denote the \textit{nodal} coordinates in physical space and time and $\xi$ and $\tau$ 
are the reference coordinates. A sketch of this isoparametric mapping is depicted in Fig. \ref{fig.isoparam}. 
Since the time coordinates are the same for each spatial node, one gets the following simple mapping for the time coordinate: 
\begin{equation}
\label{eqn.time.map} 
 t = t^n + \Delta t \tau,  
\end{equation} 
which reduces the Jacobian of the space--time mapping $(\xi,\tau) \to (x,t)$ given by \eqref{eqn.iso.map} to  
\begin{equation}
 J = \left( \begin{array}{cc} x_\xi & x_\tau \\ t_\xi & t_\tau \end{array} \right) = \left( \begin{array}{cc} x_\xi & x_\tau \\ 0 & \Delta t \end{array} \right)   
\label{eqn.jac} 
\end{equation}  
and its inverse is given by 
% Jacobi determinant 
%\begin{equation}
%|J| = \Delta t x_\xi = \Delta t \frac{\partial \theta_m}{\partial \xi} \hat x_m. 
%\end{equation}
\begin{equation}  
J^{-1} = \left( \begin{array}{cc} \xi_x & \xi_t \\ \tau_x & \tau_t  \end{array} \right)  
       = \left( \begin{array}{cc} \frac{1}{x_\xi} & - \frac{1}{\Delta t} \frac{x_\tau}{x_\xi} \\ 0 & \frac{1}{\Delta t} \end{array} \right). 
\label{eqn.invjac} 
\end{equation}      
 
Now, the derivatives of the PDE \eqref{eqn.pde} are transformed to derivatives with respect to the reference element using the chain rule, i.e. we get 
\begin{equation}
  \frac{\d \Q}{\d \tau} \tau_t + \frac{\d \Q}{\d \xi} \xi_t + \frac{\d \F}{\d \xi} \xi_x + \frac{\d \F}{\d \tau} \tau_x = \S(\Q),   
\label{eqn.pde.ref.pre}  
\end{equation} 
and with the inverse of the Jacobian of the mapping \eqref{eqn.invjac} one obtains the PDE rewritten in reference coordinates: 
\begin{equation}
  \frac{\d \Q}{\d \tau} + \frac{\Delta t}{x_\xi} \frac{\d \F}{\d \xi}  - \frac{x_\tau}{x_\xi} \frac{\d \Q}{\d \xi}  = \Delta t \S(\Q).   
\label{eqn.pde.ref}  
\end{equation} 
To simplify the notation, we introduce the following operators on the space-time reference element $T_E$: 
\begin{equation}
\label{eqn.op.ref} 
  \left[f,g\right]^\tau  = \int \limits_0^1 f(\xi,\tau) g(\xi,\tau) \, d\xi, \qquad  \textnormal{ and } \qquad 
  \left<f,g\right>       = \int \limits_{0^+}^1 \int \limits_0^1 f(\xi,\tau) g(\xi,\tau) \, d\xi \, d\tau.  
\end{equation}  
For isoparametric elements, the \textit{discrete} solution of PDE \eqref{eqn.pde.ref} is approximated using the same space-time basis functions $\theta_m$ that have also
been used for the mapping \eqref{eqn.iso.map}, i.e. 
\begin{equation}
\label{eqn.qh} 
\q_h = \q_h(x,t) = \theta_m(x,t) \hat \q_m.  
\end{equation} 
Eqn. \eqref{eqn.pde.ref} is now multiplied with space--time test functions $\theta_k$ (the same as the basis functions for the discrete solution and the mapping), and is 
then integrated over the space-time reference element $T_E = [0;1]^2$.  
One obtains 
\begin{equation}
\label{eqn.weak1} 
\left< \theta_k, \frac{\partial}{\partial \tau} \q_h \right> + \left[ \theta_k, \q_h - \w_h \right]^{0} + 
\left< \theta_k, \frac{\Delta t}{x_\xi}  \frac{\partial}{\partial \xi} \F_h(\q_h) - \frac{x_\tau}{x_\xi} \frac{\d \q_h}{\d \xi} \right> = 
\Delta t \left< \theta_k, \S_h(\q_h) \right>.  
\end{equation}
Here, the initial condition given by the WENO reconstruction polynomials $w_h(x,t^n)$ at time $t^n$ has been introduced in a \textit{weak form} via the jump term 
$\left[ \theta_k, \q_h - \w_h \right]^0$. 
In other words, the first term of \eqref{eqn.weak1} is the integral over the smooth part of the solution while the second term takes into account the jump of the discrete 
solution from $t^n$ to $t^{n,+}$. Since we use a \textit{nodal} basis, the interpolation polynomials for the flux and source terms are given by 
\begin{equation}
\label{eqn.FhSh} 
\F_h = \F_h(\q_h)  = \theta_m(\xi,\tau) \hat \F_m, \quad  \textnormal{ and }  \quad  \S_h = \S_h(\q_h)  = \theta_m(\xi,\tau) \hat \S_m,     
\end{equation} 
with 
\begin{equation}
\label{eqn.nodalapprox} 
 \hat \F_m = \F(\hat \q_m)  \quad \textnormal{ and }  \quad  \hat \S_m = \S(\hat \q_m). 
\end{equation} 

Using the definitions for the WENO reconstruction polynomial \eqref{eqn.weno} and the definitions for the discrete space--time solutions \eqref{eqn.qh}, \eqref{eqn.FhSh} 
and \eqref{eqn.nodalapprox} we obtain the following element--local nonlinear algebraic equation system:  
\begin{equation}
\label{eqn.weak3} 
 K^1_{km} \hat \q_m + K^{\xi_x}_{km} \hat \F_m - K^{\xi_t}_{km} \hat \q_m = F^0_{km} \hat \w_m^n + \Delta t M_{km} \hat \S_m,  
\end{equation} 
with the definitions of the following matrices:  
\begin{equation}
  K^1_{km} = \left( \left< \theta_k,  \frac{\partial}{\partial \tau} \theta_m \right> + \left[ \theta_k, \theta_m \right]^0 \right),  
\end{equation}  
\begin{equation}
 K^{\xi_x}_{km} = \left< \theta_k,  \frac{\Delta t}{x_\xi} \frac{\partial \theta_m}{\partial \xi} \right>, \qquad 
 K^{\xi_t}_{km} = \left< \theta_k,  \frac{x_\tau}{x_\xi}   \frac{\partial \theta_m}{\partial \xi} \right>, 
\end{equation}  
and 
\begin{equation}
F^0_{km} = \left[\theta_k, \psi_m \right]^0, \qquad 
M_{km} = \left< \theta_k, \theta_m  \right>. 
\end{equation} 
Due to the nodal approach on the tensor--product Gauss--Legendre nodes, the $\xi$ and $\tau$ directions decouple and the above matrices can be easily evaluated 
dimension by dimension using one-dimensional Gaussian quadrature. 
The system \eqref{eqn.weak3} is solved using an iterative method similar to the one proposed in \cite{DumbserZanotti,HidalgoDumbser}: 
\begin{equation}
\label{eqn.weak.iter} 
 K^1_{km} \hat \q_m^{l+1} + K^{\xi_x}_{km} \hat \F_m^{l} - K^{\xi_t}_{km} \hat \q_m^{l} = F^0_{km} \hat \w_m^n + \Delta t M_{km} \hat \S_m^{l+1},  
\end{equation} 
where the stiff algebraic source term is taken implicitly
(see \cite{DumbserZanotti} for details). A particularly efficient strategy for obtaining an initial guess for 
$\hat \q^0_m$ can be found in \cite{HidalgoDumbser}. 
The equation that determines the location of the spatial coordinates $\hat x_m$ of the space--time element is
\begin{equation}
\label{eqn.mesh.motion} 
  \frac{dx}{dt} = V(x,t),  
\end{equation} 
where $V(x,t)$ is the local mesh velocity. In the fully Eulerian case one has $V=0$, while in the pure Lagrangian case, $V(x,t)$ is the local 
fluid velocity. For the local mesh velocity we use the nodal ansatz 
\begin{equation}
  V_h = V_h(x,t) = \theta_m(x,t) \hat v_m, \qquad \textnormal{ with } \qquad \hat v_m = V(\hat x_m,\hat t_m). 
\end{equation}
The initial distribution of the spatial Gauss--Legendre quadrature points at time $t^n$ is given by
\begin{equation}
  \hat x_{0,m} = x_{i-\halb}^n + \Delta x_i^n \zeta_m, 
\end{equation} 
and the spatial Lagrange interpolation polynomials passing through these points are denoted by $\phi_m$. 
Formally, a discrete version of the ODE \eqref{eqn.mesh.motion} can then be obtained using again the local space--time DG method, 
see \cite{ADERNSE}: 
\begin{equation}
\label{eqn.weak.x} 
  \left( \left< \theta_k, \frac{\partial}{\partial \tau} \theta_m \right> + [\theta_k, \theta_m]^0 \right) \hat x^{l+1}_m = [\theta_k, \phi_m]^0 \hat x_{0,m} + 
   \Delta t \left< \theta_k, \theta_m \right> \hat v^l_m.  
\end{equation} 
The weak formulation for the spatial coordinates \eqref{eqn.weak.x} is iterated \textit{together} with the weak formulation for the solution \eqref{eqn.weak.iter} until 
convergence is reached. The temporal coordinates $\hat t_m$ are always \textit{fixed} and are given by the Gauss--Legendre points $\zeta_j$ and the relation \eqref{eqn.time.map}.  
The space--time polynomials $\q_h(x,t)$ are computed for each element in the computational domain and are then used as arguments for the numerical flux in Eqn. \eqref{eqn.numflux}
and the numerical source term $\S_i$ is computed using $\q_h(x,t)$ as follows:
\begin{equation}
  \S_i = \frac{1}{\Delta x_i^{n} \, \Delta t} \int \limits_{t^n}^{t^{n+1}} \int \limits_{x_{i-\halb}(t)}^{x_{i+\halb}(t)} \S\left(\q_h(x,t)\right) \ dx dt.     
\end{equation}  
This completes the description of the high order Lagrangian one--step finite volume algorithm \eqref{eqn.fv}. 

%\clearpage 
%%%DEBUG
%\begin{equation}
%\label{eqn.weak2} 
%\left< |J| \theta_k,  \left( \xi_t \frac{\partial}{\partial \xi} + \tau_t \frac{\partial}{\partial \tau}   \right) \q_h  \right> + \left[ x_\xi \theta_k, \q_h - \w_h \right]^0 + 
%\left< |J| \theta_k,  \xi_x \frac{\partial \F_h}{\partial \xi} \right> = \left< |J| \theta_k, \S_h  \right>.  
%\end{equation} 
%\begin{equation} 
%\left< \theta_k, \frac{\partial}{\partial \tau} \theta_l \right> (|J| \tau_t \hat q)_l + 
%\left< \theta_k, \frac{\partial}{\partial \xi}  \theta_l \right> (|J| \xi_t  \hat q)_l + 
%\left[ \theta_k, \theta_l \right]^0 (x_\xi \hat q)_l - \left[ \theta_k, \psi_l \right]^0 (x_\xi \hat w)_l -  
%\left< |J| \theta_k,  \xi_x \frac{\partial \F_h}{\partial \xi} \right> = \left< \theta_k, \theta_l \right> (|J| \hat \S)_l.  
%\end{equation} 
%\begin{equation} 
%\left< \theta_k, \frac{\partial}{\partial \tau} \theta_l \right> (|J| \hat q)_l + 
%\left< \theta_k, \frac{\partial}{\partial \xi}  \theta_l \right> (|J| \Delta t \xi_t  \hat q)_l + 
%\left[ \theta_k, \theta_l \right]^0 (|J| \hat q)_l - \left[ \theta_k, \psi_l \right]^0 (|J| \hat w)_l -  
%\left< \Delta t |J| \theta_k,  \xi_x \frac{\partial \F_h}{\partial \xi} \right> = \Delta t \left< \theta_k, \theta_l \right> (|J| \hat \S)_l.  
%\end{equation} 
%%%DEBUG
%\clearpage 

\section{Test problems} 
\label{sec.test} 

In this section we show some computational test problems to illustrate the performance of the scheme in the case of the compressible Euler equations (non--stiff), 
the resistive relativistic MHD equations (stiff) and for the relativistic radiation hydrodynamics equations (moderately stiff).  

\subsection{Compressible Euler equations} 
\label{sec.euler} 

The Euler equations of compressible gas dynamics read 
\begin{equation}
\label{eqn.euler} 
  \frac{\partial}{\partial t} \left( \begin{array}{c} \rho \\ \rho u \\ \rho E \end{array} \right) + 
  \frac{\partial}{\partial x} \left( \begin{array}{c} \rho u \\ \rho u^2 + p \\ u( \rho E + p) \end{array} \right) = \S(x,t),     
\end{equation}
with the fluid density $\rho$, the velocity $u$, the total energy density $\rho E$, a vector of source terms $\S$ and the fluid pressure $p$, 
given in terms of the conserved quantities by the equation of state of an ideal gas as 
\begin{equation}
\label{eqn.eos} 
p = (\gamma-1)(\rho E - \halb \rho u^2),  
\end{equation}
with the ratio of specific heats $\gamma$. In this section, we define the local mesh velocity as the local fluid velocity, i.e. we choose 
\begin{equation}
 V = u. 
\end{equation} 

\subsubsection{Numerical convergence results} 
\label{sec.conv.Euler}

In order to assess the accuracy of the method presented in section \ref{sec.method} we carry out several simulations of a test problem with smooth solution and exact solution on a series of successively refined meshes. For this purpose, we use the so--called manufactured solution method, which means prescribing the exact solution a priori, inserting the solution in PDE \eqref{eqn.euler} and putting all terms that do not cancel into the source term $\S(x,t)$. For our particular test problem we choose the primitive variables of the exact solution as 
\begin{equation}
\label{eqn.euler.exact} 
  \rho_e(x,t) = 1 + \frac{1}{2} \sin(2 \pi x) \cos(2 \pi t), \qquad u_e(x,t) = \sin(2 \pi x) \cos(2 \pi t), \qquad p_e(x,t) = 1.  
\end{equation} 
From there, the vector of conserved variables can be computed as $\Q_e = (\rho_e, \rho_e u_e, p_e/(\gamma-1) + 1/2 \rho_e u_e^2)^T$. 
The ratio of specific heats is chosen as $\gamma=1.4$. 
Inserting \eqref{eqn.euler.exact} into \eqref{eqn.euler} yields the source term $\S(x,t)$. The initial computational domain is 
$\Omega = [-2;2]$ and is discretized by an initially uniform mesh of $N_G$ control volumes. The boundary conditions are periodic. 
Simulations are carried out with third to eighth order Lagrangian one--step WENO finite volume schemes using the Osher--type flux \eqref{eqn.osher} 
for one period until the final time $t=1.0$. The Courant number is set to $CFL=0.9$. For each mesh the corresponding error in 
$L_2$ norm is computed as 
\begin{equation}
  \epsilon_{L_2} = \sqrt{ \int \limits_{\Omega(t_e)} \left( \Q_e(x,t_e) - \w_h(x,t_e) \right)^2 dx },  
\end{equation} 
and the resulting numerical convergence rates are listed in Table \ref{tab.conv1}. From the presented results we can conclude that the
designed order of accuracy of the scheme is reached.   

\begin{table}  
\caption{Numerical convergence results for the compressible Euler equations using the third to eighth order version of the Lagrangian one--step 
WENO finite volume schemes presented in this article. The error norms refer to the variable $\rho$ (density) at time $t=1.0$.} 
\begin{center} 
\renewcommand{\arraystretch}{1.0}
\begin{tabular}{ccccccccc} 
\hline
  $N_G$ & $\epsilon_{L_2}$ & $\mathcal{O}(L_2)$ & $N_G$ & $\epsilon_{L_2}$ & $\mathcal{O}(L_2)$ & $N_G$ & $\epsilon_{L_2}$ & $\mathcal{O}(L_2)$ \\ 
\hline
  \multicolumn{3}{c}{$\mathcal{O}3$} & \multicolumn{3}{c}{$\mathcal{O}4$}  & \multicolumn{3}{c}{$\mathcal{O}5$} \\
\hline
100  & 1.9526E-02 &     & 100  & 5.7994E-03 &     &  50  & 2.1944E-02 &      \\ 
200  & 3.0021E-03 & 2.7 & 200  & 1.2551E-04 & 5.5 & 100  & 1.5756E-03 & 3.8  \\ 
400  & 4.0927E-04 & 2.9 & 400  & 4.9135E-06 & 4.7 & 200  & 8.2557E-05 & 4.3  \\ 
800  & 5.7539E-05 & 2.8 & 800  & 3.1365E-07 & 4.0 & 400  & 3.3144E-06 & 4.6  \\ 
\hline 
  \multicolumn{3}{c}{$\mathcal{O}6$} & \multicolumn{3}{c}{$\mathcal{O}7$}  & \multicolumn{3}{c}{$\mathcal{O}8$} \\
\hline
 50  & 1.6783E-02 &     &  50  & 9.2252E-03 &     &  50  & 6.9643E-03 &      \\ 
100  & 6.5205E-04 & 4.7 & 100  & 2.7828E-04 & 5.1 & 100  & 1.6253E-04 & 5.4  \\ 
200  & 7.4380E-06 & 6.5 & 200  & 5.6255E-06 & 5.6 & 200  & 6.9410E-07 & 7.9  \\ 
400  & 9.1756E-08 & 6.3 & 400  & 6.8116E-08 & 6.4 & 300  & 2.7163E-08 & 8.0  \\ 
\hline 
\end{tabular} 
\end{center}
\label{tab.conv1}
\end{table} 
	
\subsubsection{Shock tube problems} 
\label{sec.shock.Euler}

\begin{figure}[!htbp]
\begin{center}
\begin{tabular}{lr}
\includegraphics[width=0.45\textwidth]{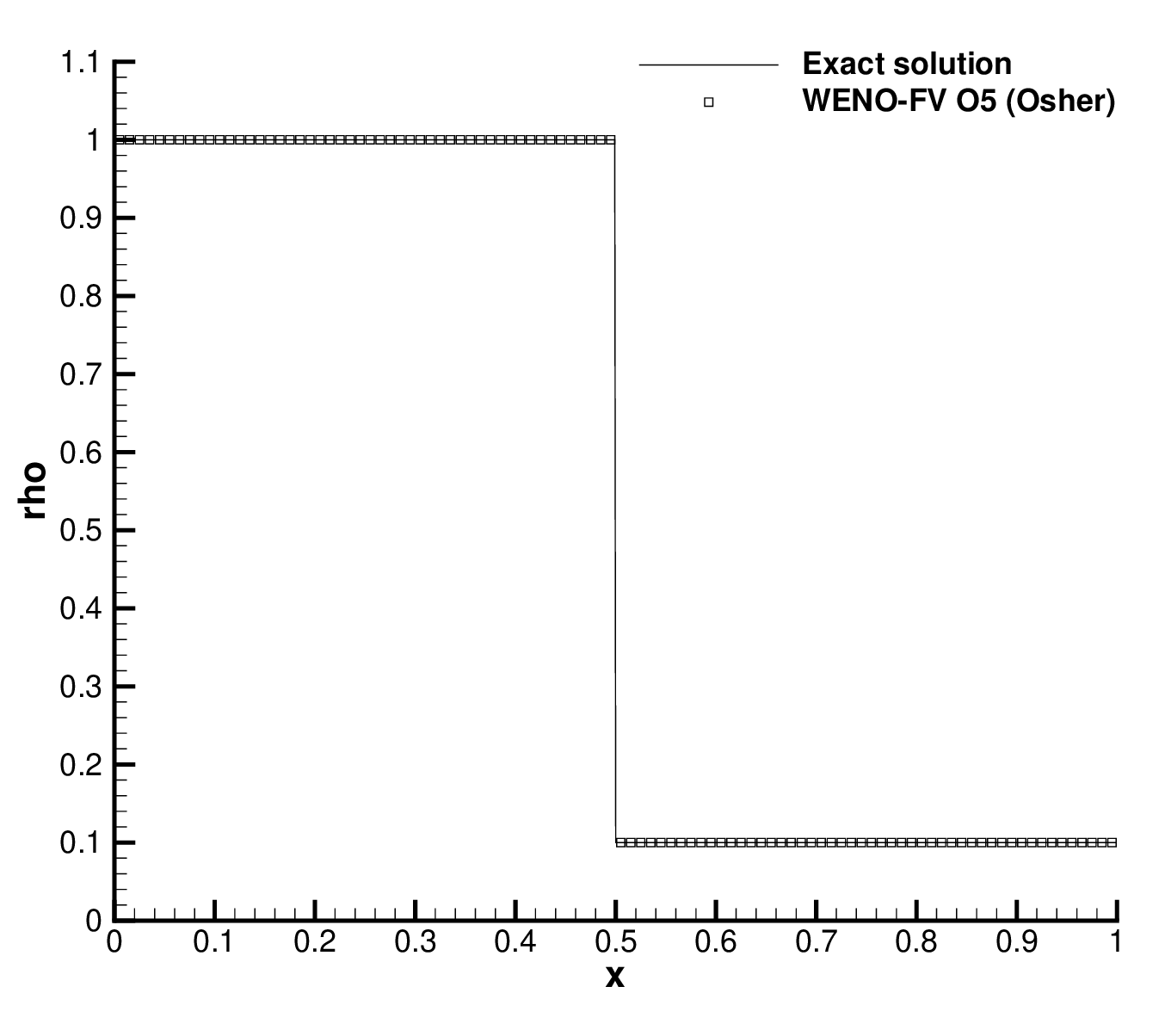}   & 
\includegraphics[width=0.45\textwidth]{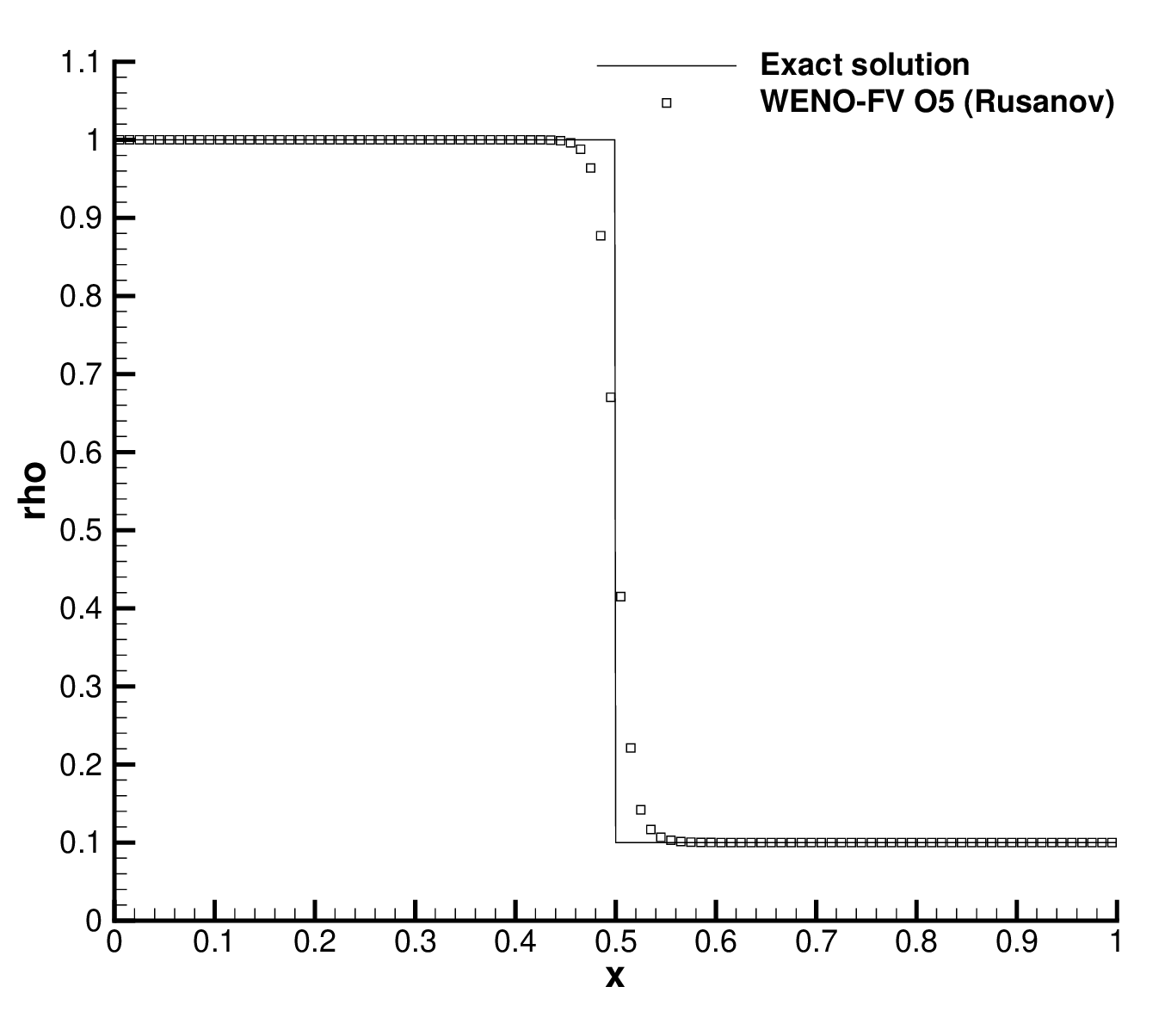} 
\end{tabular}
\caption{Exact and numerical solution obtained with third order Lagrangian one--step WENO finite volume schemes for RP0 
(isolated moving contact wave) at $t=0.5$ using 100 cells. Left: Osher--type flux \eqref{eqn.osher}. Right: Rusanov--type flux \eqref{eqn.rusanov}.  }
\label{fig.rp0}
\end{center}
\end{figure}

In this section we solve a set of several shock--tube problems given in \cite{toro-book} and \cite{chengshu1,chengshu2}. 
The initial conditions of the Riemann problems are 
\begin{equation}
  {\bf Q}(x,0) = \begin{cases} {\bf Q}_L & \IF x<x_d, \\ {\bf Q}_R & \IF x \geq x_d, \end{cases}  
\end{equation} 
where the initial states left and right are summarized in Table \ref{tab.ic.euler}. The initial computational 
domain is $\Omega = [x_L; x_R]$ and is discretized with 200 equidistant cells, apart from RP5 (Leblanc shock tube), for
which 2000 cells are used, and RP0, for which only 100 cells are used. Simulations have been carried out with the fifth order 
version of our Lagrangian one--step WENO finite volume schemes. % The Courant number has been set to CFL=0.5 in all cases. 
The first problem (RP0) is the advection of an isolated moving contact wave with constant pressure and velocity. Any
Riemann solver that resolves exactly stationary contact waves in the Eulerian case should preserve exactly isolated moving 
contact waves in the Lagrangian case. From Fig. \ref{fig.rp0} it becomes evident that the Osher-type flux \eqref{eqn.osher} 
solves the problem \textit{exactly}, without any intermediate points in the contact wave, whereas the Rusanov flux \eqref{eqn.rusanov} 
adds significant numerical diffusion to the problem, as expected. For this problem, the Osher-type flux leads to a pure Lagrangian 
scheme, where the mass in each moving control volume remains constant.  

\begin{figure}[!htbp]
\begin{center}
\begin{tabular}{lr}
\includegraphics[width=0.45\textwidth]{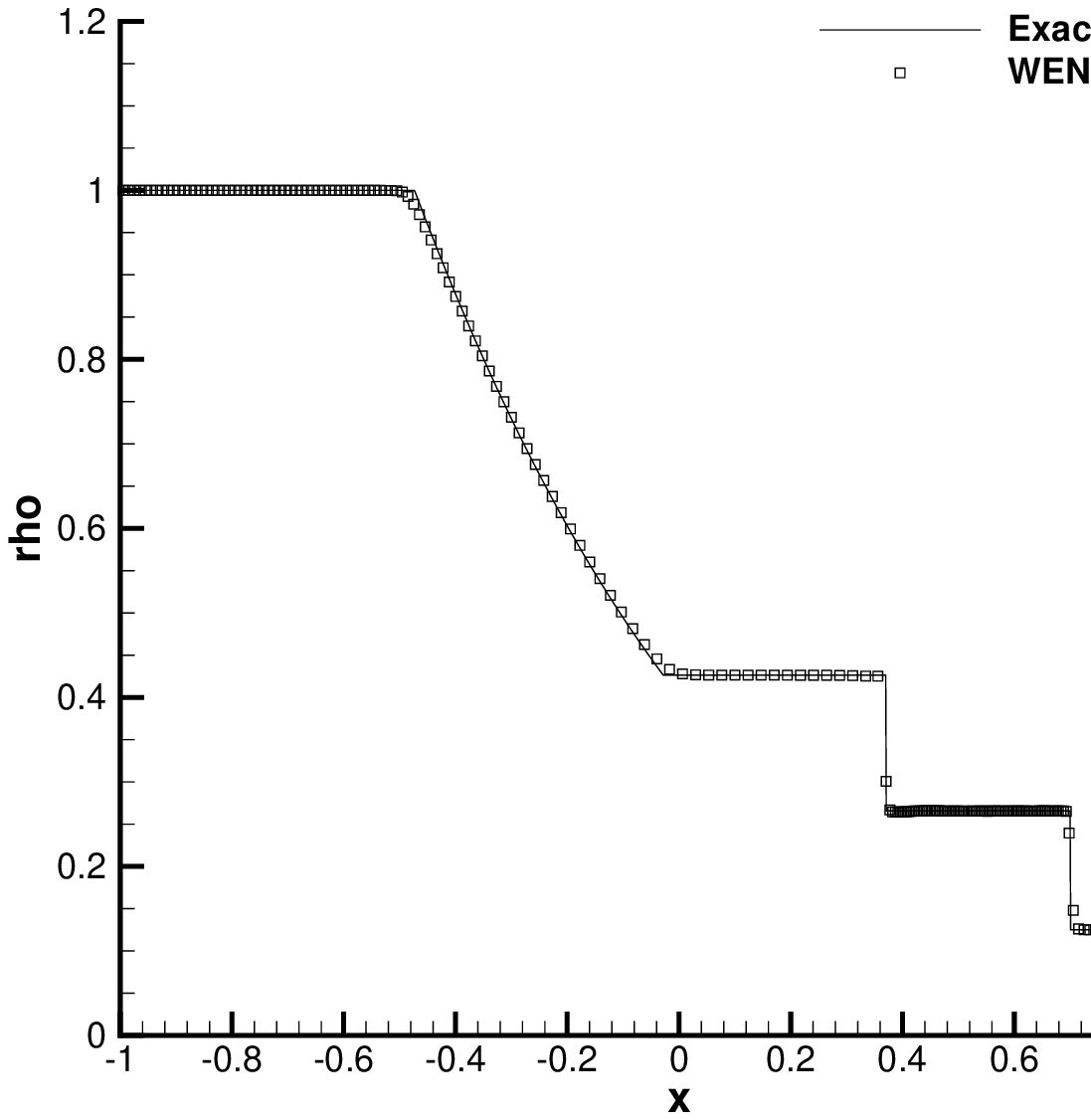} &
\includegraphics[width=0.45\textwidth]{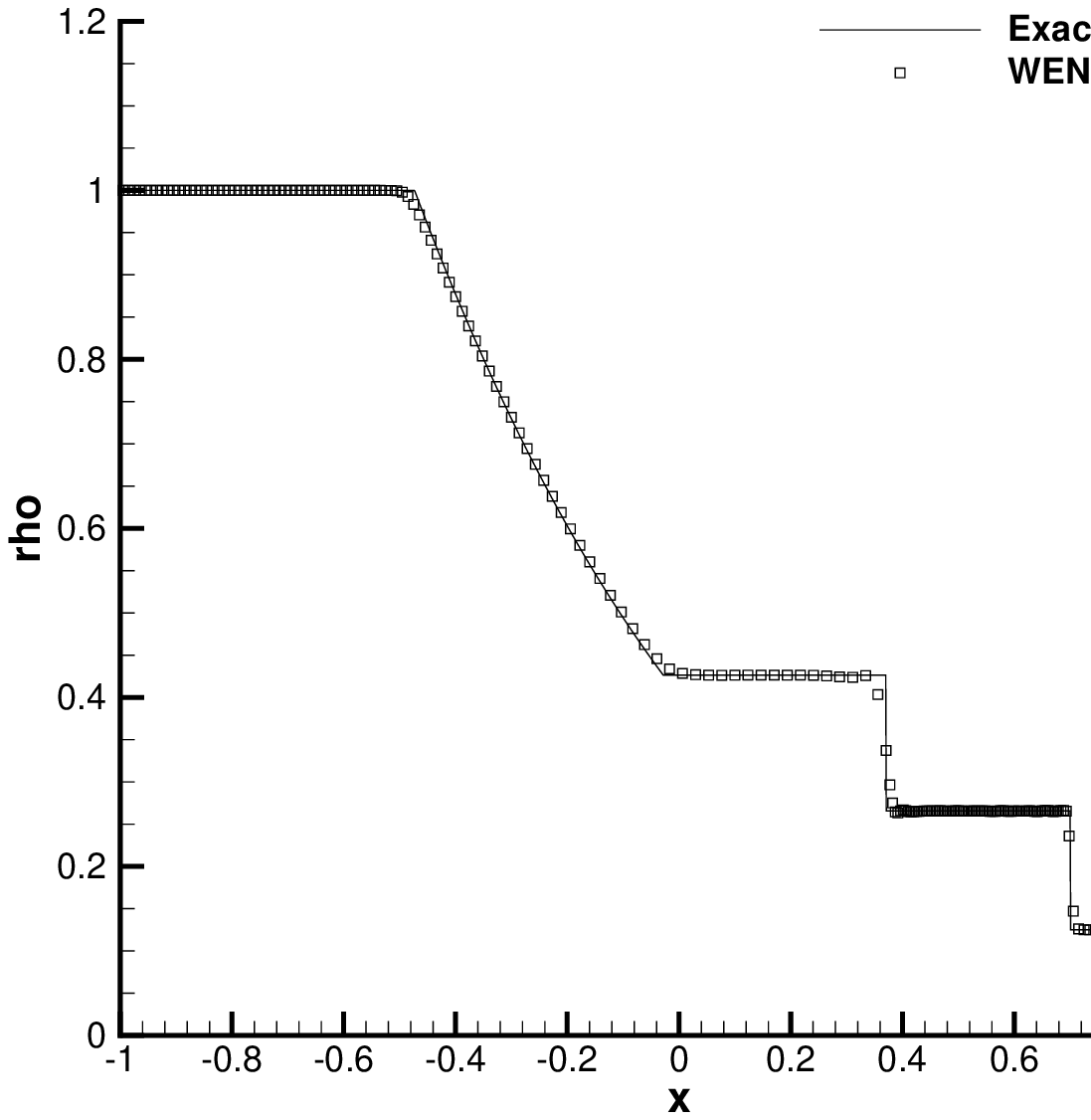} \\ 
\includegraphics[width=0.45\textwidth]{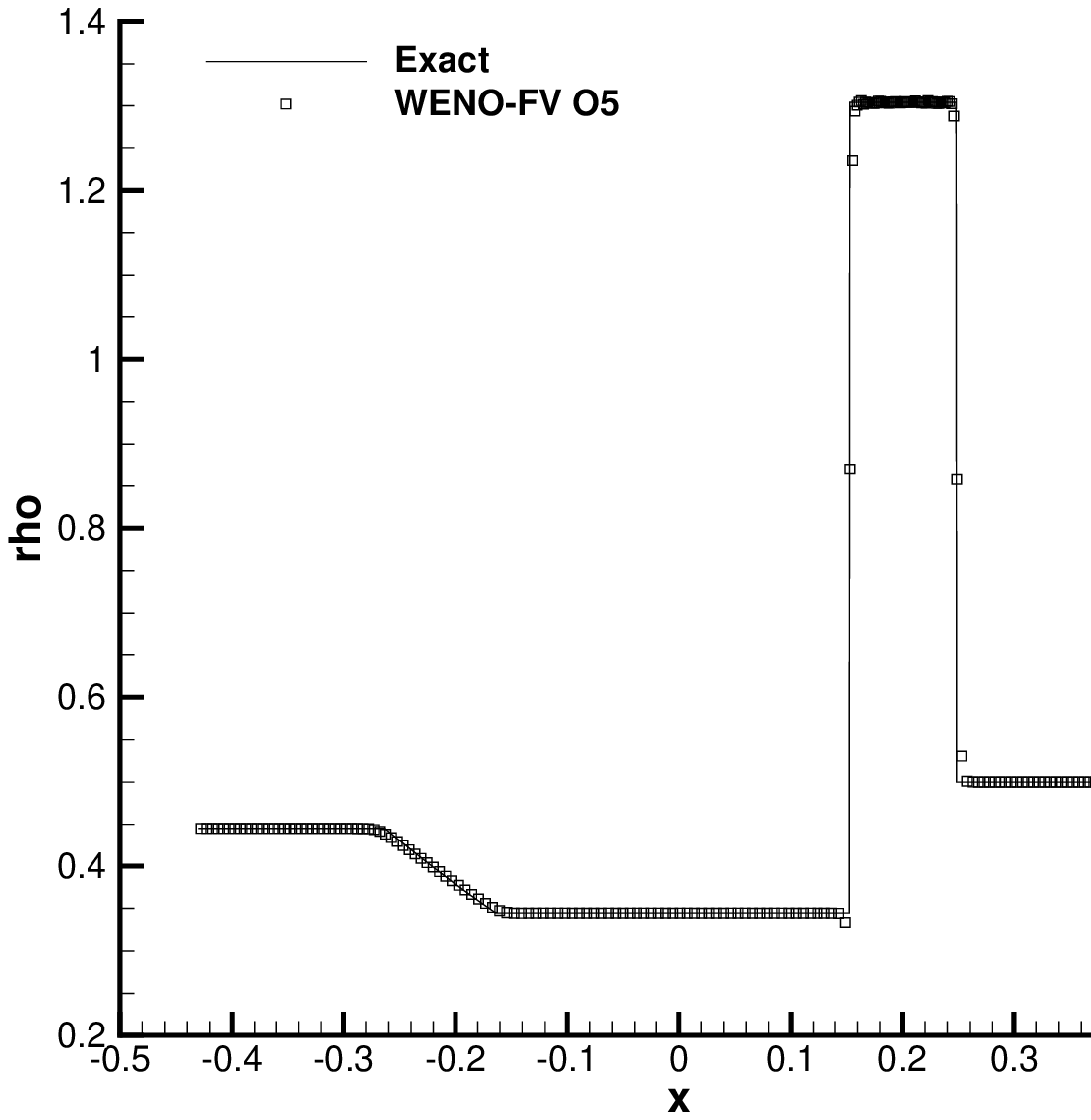} &
\includegraphics[width=0.45\textwidth]{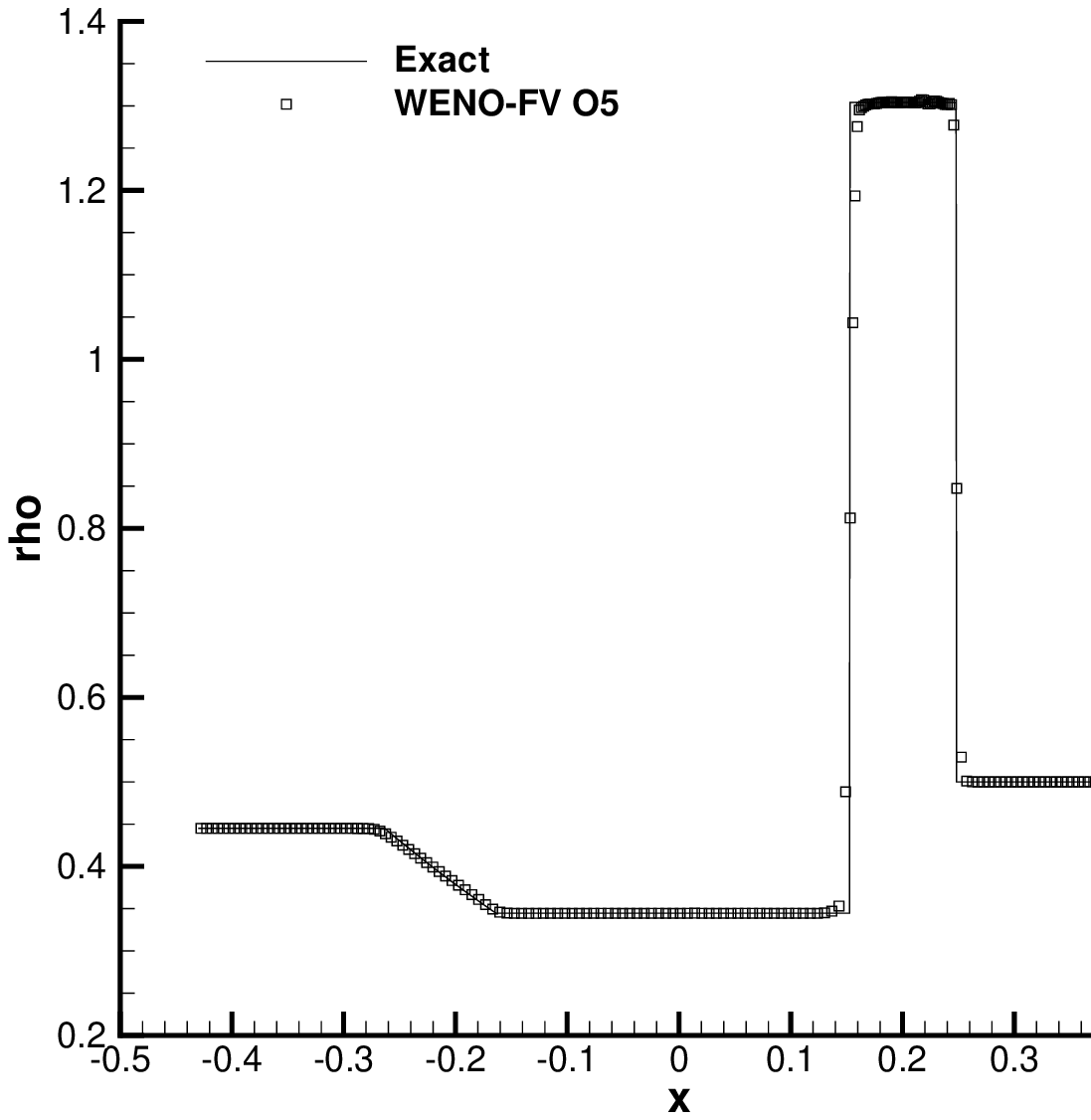} \\  
\end{tabular}
\caption{Exact and numerical solution obtained with fifth order Lagrangian one--step WENO finite volume schemes for RP1 
(Sod problem) at $t=0.4$ (top) and RP2 (Lax problem) at $t=0.1$ (bottom). 
Left: Osher--type flux \eqref{eqn.osher}. Right: Rusanov--type flux \eqref{eqn.rusanov}. }
\label{fig.sodlax}
\end{center}
\end{figure}

\begin{figure}[!htbp]
\begin{center}
\begin{tabular}{lr}
\includegraphics[width=0.45\textwidth]{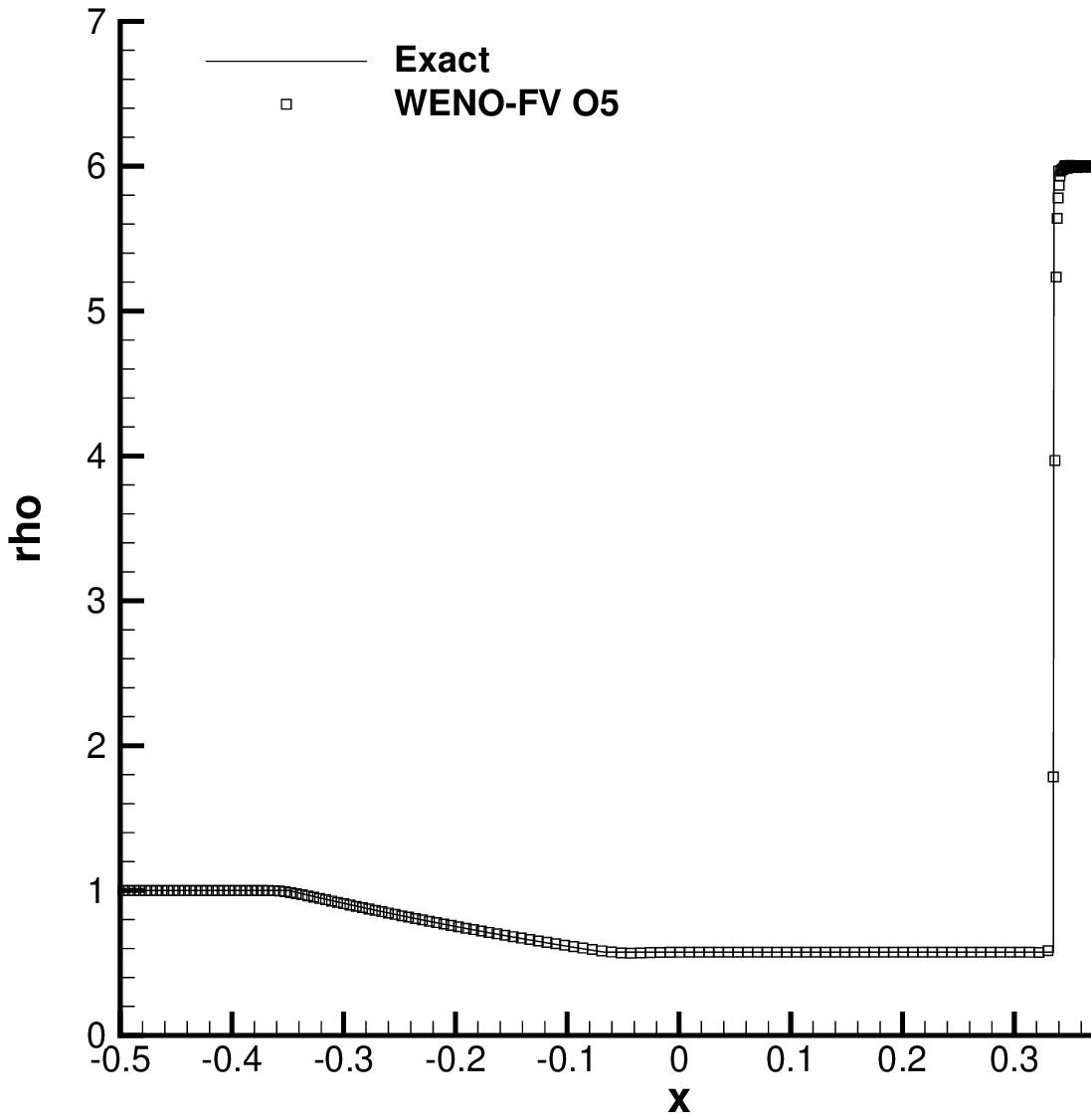} &
\includegraphics[width=0.45\textwidth]{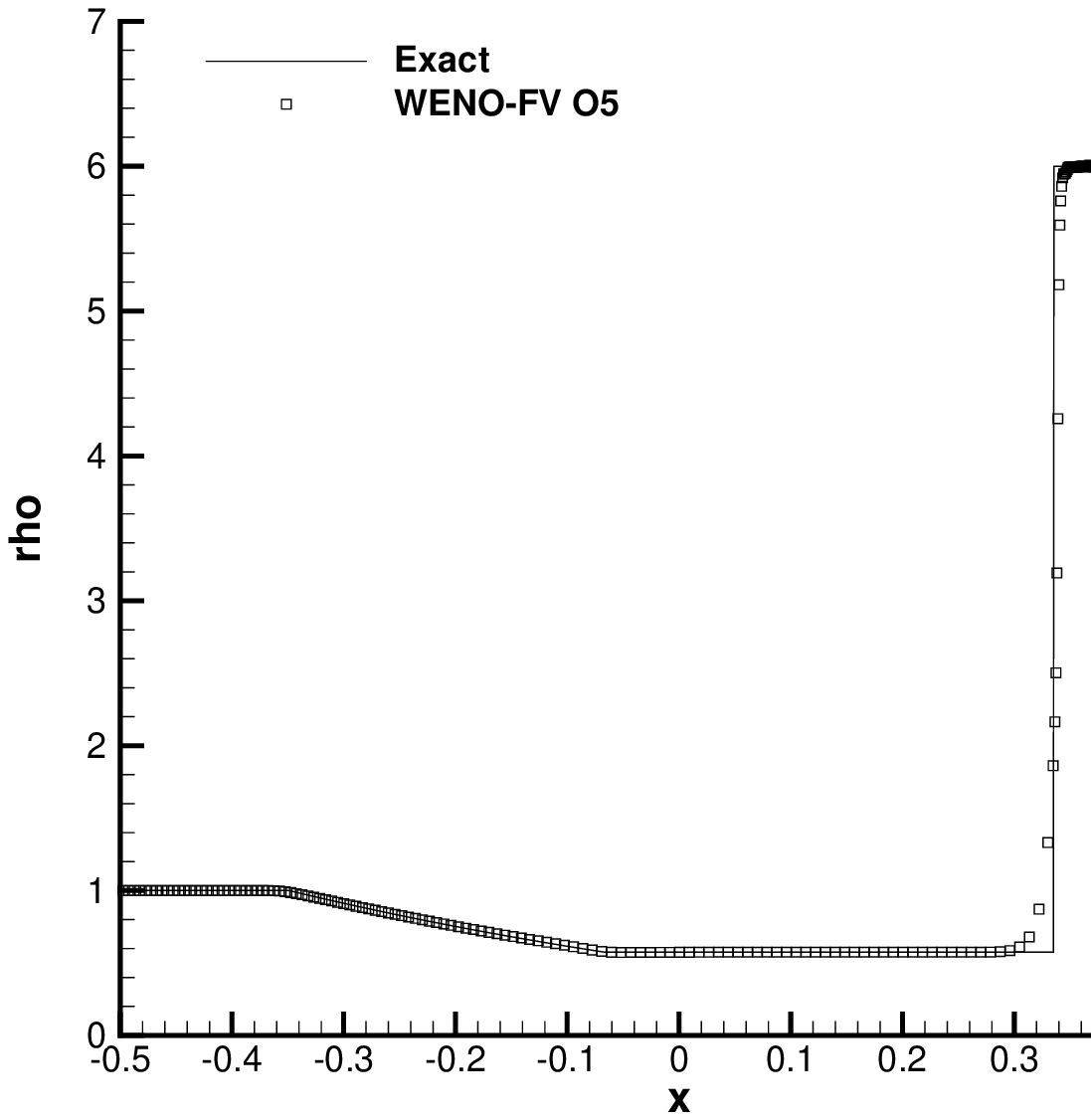} \\ 
\includegraphics[width=0.45\textwidth]{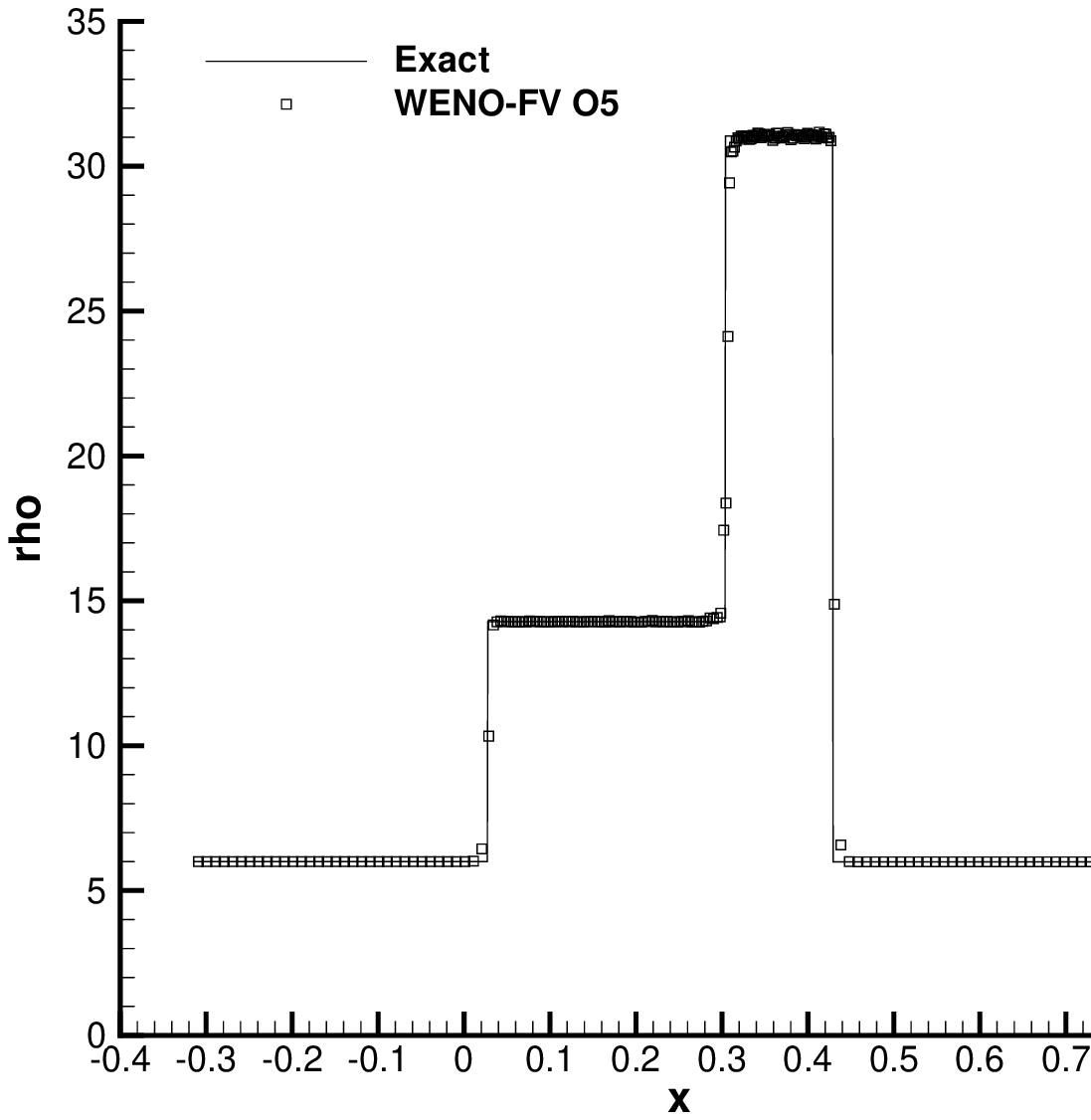} &
\includegraphics[width=0.45\textwidth]{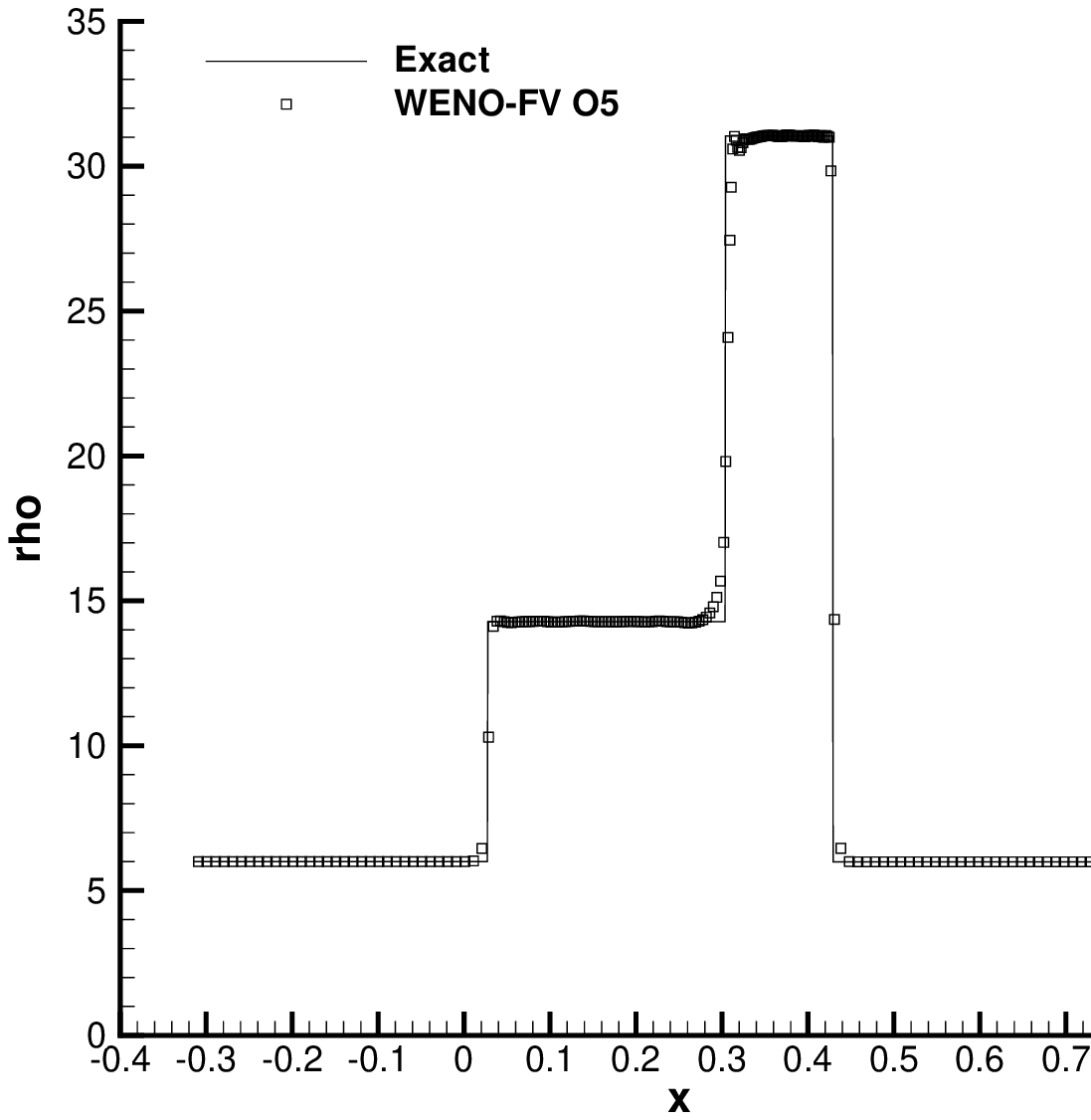} \\  
\end{tabular}
\caption{Exact and numerical solution obtained with fifth order Lagrangian one--step WENO finite volume schemes for RP3 
at $t=0.012$ (top) and RP4 at $t=0.035$ (bottom). 
Left: Osher--type flux \eqref{eqn.osher}. Right: Rusanov--type flux \eqref{eqn.rusanov}. }
\label{fig.toro34}
\end{center}
\end{figure}

\begin{figure}[!htbp]
\begin{center}
\begin{tabular}{lr}
\includegraphics[width=0.45\textwidth]{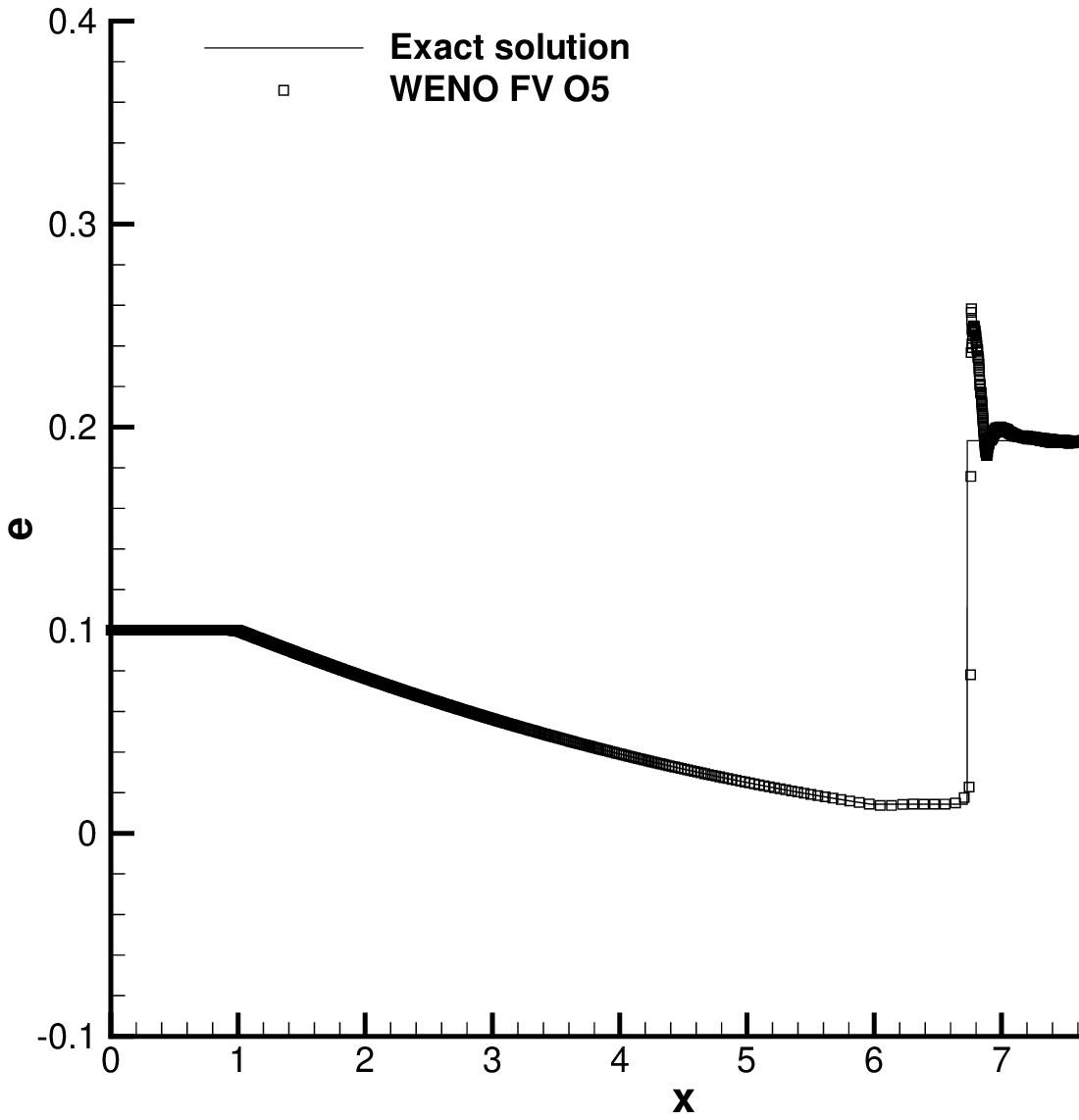} &
\includegraphics[width=0.45\textwidth]{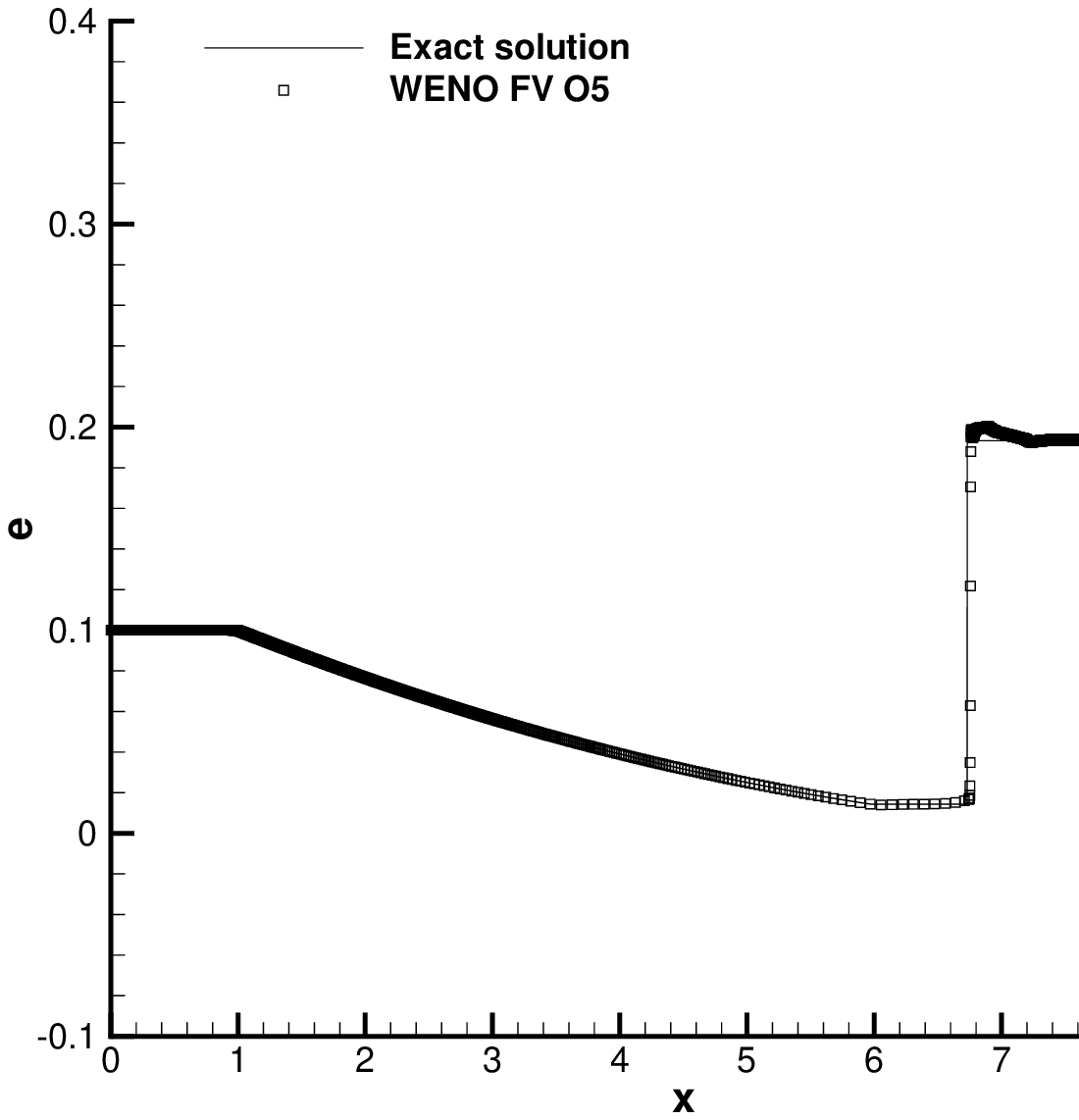} \\ 
\end{tabular}
\caption{Exact and numerical solution obtained with fifth order Lagrangian one--step WENO finite volume schemes for RP5 
(Leblanc problem) at $t=6.0$. Left: Osher--type flux \eqref{eqn.osher}. Right: Rusanov--type flux \eqref{eqn.rusanov}. }
\label{fig.leblanc}
\end{center}
\end{figure}

In Figs. \ref{fig.sodlax} - \ref{fig.leblanc} we show the exact solution together with the computational results obtained 
with Osher--type flux \eqref{eqn.osher} and the Rusanov--type flux \eqref{eqn.rusanov} for the other shock tube problems 
RP1-RP5. Overall, a very good agreement is noted between the numerical solution and the exact solution. The Osher flux resolves 
the contact wave very well in general. However, some intermediate points are generated, since in the initial phase of the Riemann 
problem waves of different nature (shock and rarefaction) overlap with the contact wave, thus leading to some amount of numerical 
diffusion in the contact wave. 
These results are as expected, since in the present paper an ALE-type approach has been used, which does in general \textit{not} 
impose constant mass in each control volume, in contrast to purely Lagrangian schemes as the one presented, for example, in 
\cite{Carre2009}.  
For the Leblanc problem, it can be easily noted that the Rusanov flux is much more robust for this problem due to its larger 
numerical diffusion compared to the Osher-type scheme. However, the results presented here are similar to the ones presented 
in \cite{chengshu1,chengshu2}. 

\begin{table}[!t]
 \caption{Initial states left and right for density $\rho$, velocity $u$ and pressure $p$ for the Riemann problems 
 solved for the compressible Euler equations. The initial position of the discontinuity ($x_d$) and the initial computational 
 domain $\Omega=[x_L;x_R]$ are also specified. In all cases $\gamma=1.4$, apart for RP5, where $\gamma=5/3$. } 
\begin{center} 
 \begin{tabular}{cccccccccc} 
 \hline
 RP & $\rho_L$ & $u_L$ & $p_L$ & $\rho_R$ & $u_R$ & $p_R$ & $x_d$ & $x_L$ & $x_R$  \\ 
 \hline
 0 &  1.0      &  1.0       & 1.0     & 0.1        &  1.0        & 1.0      &   0.0    & -0.5 & 0.5 \\
 1 &  1.0      &  0.0       & 1.0     & 0.125      &  0.0        & 0.1      &   0.0    & -1.0 & 1.0 \\
 2 &  0.445    &  0.698     & 3.528   & 0.5        &  0.0        & 0.571    &   0.0    & -0.5 & 0.5 \\ 
 3 &  1.0      &  0.0       & 1000    & 1.0        &  0.0        & 0.01     &   0.1    & -0.5 & 0.5 \\
 4 &  5.99924  &  19.5975   & 460.894 & 5.99242    & -6.19633    & 46.095   &   0.0    & -1.0 & 1.0 \\ 
 5 &  1.0      &  0.0       & $0.1(\gamma-1)$      & $10^{-3}$   & 0.0      & $10^{-10}(\gamma-1)$ &  3.0    &  0.0 & 9.0 \\
 \hline
 \end{tabular}
\end{center} 
 \label{tab.ic.euler}
\end{table}

\subsection{Resistive relativistic MHD equations (RRMHD)} 
\label{sec.rrmhd} 

The resistive relativistic MHD (RRMHD) equations constitute a hyperbolic system of balance laws which 
has a source term that may become stiff. The equations include five equations for the fluid motion 
(conservation of mass, momentum and energy), plus six equations for the evolution of the electric and of the 
magnetic field (Maxwell equations). Furthermore, two additional equations are needed to maintain the constraints 
on the divergence of the electric and of the magnetic field. In this paper, we use the hyperbolic divergence cleaning 
approach according to Dedner et al. \cite{Dedneretal}. The last equation expresses the conservation of the total charge. 
In Cartesian coordinates, using the abbreviations $\partial_t = \frac{\partial}{\partial t}$ and $\partial_i = \frac{\partial}{\partial x_i}$, 
the resistive relativistic MHD equations can be written
as follows \cite{DumbserZanotti}:  

\bea
\label{fluid1}
&& \partial_t D + \partial_i (D v^i)=0, \\
\label{fluid2-4}
&&\partial_t S_j + \partial_i Z_{j}^i=0, \\
\label{fluid5}
&&\partial_t \tau + \partial_i S^i=0, \\
\label{electric6-8}
&&\partial_t E^i - \epsilon^{ijk}\partial_j B_k + \partial_i
\Psi = -J^i, \\
&&\partial_t B^i + \epsilon^{ijk}\partial_j E_k + \partial_i
\Phi = 0, \\
&&\partial_t \Psi + \partial_i E^i = \rho_c - \kappa \Psi,
\\
&&\partial_t \Phi + \partial_i B^i = - \kappa \Phi, \\
\label{charge}
&&\partial_t \rho_c + \partial_i J^i = 0, 
\eea
where the conservative variables of the fluid are
\bea
\label{eq:D}
D&=&\rho W, \\
S^i & = & \omega W^2 v^i + \epsilon^{ijk}E_jB_k, 
\label{eq:S} \\
\tau & = & \omega W^2 - p +
\textstyle{\frac{1}{2}}(E^2+B^2) \ ,
\label{eq:U}
\eea
expressing, respectively, the relativistic mass density, the
momentum density and the total energy density.
The spatial tensor $Z^i_j$ in (\ref{fluid2-4}), representing the momentum 
flux density, is 
\begin{equation} 
\label{eq:W} 
Z^i_j  =  \omega W^2 v^i\,v_j -E^i\,E_j-B^i\,B_j+\left[p+\textstyle{\frac{1}{2}}(E^2+B^2)\right]\,\delta^i_j,
\end{equation} 
where $\delta^i_j$ is the Kronecker delta, while $W=1/\sqrt{1-v^2}$ is the Lorentz factor of the fluid.
In this paper we have assumed the equation of
state of an ideal gas, namely 
\be
\label{eos}
p=(\gamma-1) \rho\epsilon=\gamma_1(\omega-\rho), 
\ee
where $\gamma$ is the adiabatic index,
$\gamma_1=(\gamma-1)/\gamma$, $\epsilon$ is the specific internal energy and
$\omega = \rho \epsilon + \rho + p$ is the enthalpy.  
The source term $J$ appearing in \eqref{electric6-8} is
the current vector, given by  Ohm's law, for which we
assume the following expression~\cite{Komissarov2007,Palenzuela2009},   
\begin{equation} 
J^i = \rho_c v^i + \sigma W [E^i + \epsilon^{ijk} v_j B_k - (E^j v_j) v^i] \ ,
\end{equation} 
where $\rho_c$ is the charge density in the laboratory frame.

The system of equations (\ref{fluid1})-(\ref{charge}) is
written as a hyperbolic system of balance laws as in
(\ref{eqn.pde}) and it has
source terms in the three equations (\ref{electric6-8})
that are potentially stiff,  see \cite{Palenzuela2009}. 
In the stiff limit case ($\sigma \to \infty$) the resistive 
relativistic MHD equations  reduce to the ideal
relativistic MHD equations (RMHD),  for which several
test problems with exact solution are known, see
\cite{BalsaraRMHD,RMHD,GiacomazzoRezzolla}. For the
system (\ref{fluid1})-(\ref{charge}) a  family of high
order one-step schemes in \textit{Eulerian} coordinates
has been proposed in \cite{DumbserZanotti}, while in this
paper we use a \textit{Lagrangian} method.  

\subsubsection{Numerical convergence results} 
\label{sec.conv.rrmhd}

The smooth unsteady test case with exact analytical solution used here was introduced for the ideal relativistic MHD equations 
by Del Zanna et al. \cite{delZanna2007} and was solved for the first time on unstructured triangular meshes with high order
$\PNM$ schemes in \cite{Dumbser2008}. Since the resistive MHD equations tend asymptotically to the ideal ones in the stiff limit ($\sigma \to \infty$), this is an ideal test case to assess the accuracy of our scheme in the stiff limit of the governing PDE  system. \\
The test case consists of a periodic Alfv\'en wave whose initial condition at $t=0$ is chosen to be $\rho=p=1$, $B^i = B_0 \, (1, \cos\left(k x \right), \sin \left(kx \right) )^T$, 
$ v^i = -v_A/B_0 \, \cdot (0,B_y,B_z)^T$, $ E^i = -\epsilon^{ijk} v_j B_k$ and $\phi=\psi=q=0$. We furthermore use the parameters $k=2 \pi$, $\gamma=\frac{4}{3}$ and $B_0=1$, 
hence the advection speed of the Alfv\'en wave in $x$-direction is $v_A = 0.38196601125$, see \cite{delZanna2007} for a closed analytical expression for $v_A$. The computational 
domain is $\Omega = [0;1]$ with periodic boundary conditions and the final time at which we compare the exact solution with the numerical one is chosen as $t=0.5$. Since in this test
case the fluid velocity in $x$-direction is $v_x=0$, we move the mesh artificially with the fluid velocity in $y$-direction, i.e. we set $V=v_y$.  
Since this test case was constructed for the \textit{ideal} relativistic MHD equations, we have to use a very high value for the conductivity ($\sigma=10^8$) in the resistive case 
to reproduce the ideal equations asymptotically. In all our computations a \textit{constant} Courant number of $\textnormal{CFL}=0.5$ is used. \\
Table \ref{tab.conv2} shows the errors and the orders of convergence measured in the $L^2$ norm for the flow variables $v_y$ and $E_y$. 
The number $N_G$ denotes the number of grid points along the $x$-axis. We find that the nominal order of accuracy $M+1$ has been reached for
all schemes from third to sixth order of accuracy in space and time under consideration, even for the electric field $E_y$, which is one of
the variables onto which the stiff relaxation source term is acting.   
\begin{table}  
\caption{Numerical convergence results for the stiff limit ($\sigma = 10^8$) of the resistive relativistic MHD equations (RRMHD) 
using third to sixth order Lagrangian one--step WENO finite volume schemes. The error norms refer to the variable $v_y$ 
(velocity in $y$-direction) and to the \textit{relaxed} variable $E_y$ (electrical field in $y$-direction).}  
\begin{center} 
\renewcommand{\arraystretch}{1.0}
\begin{tabular}{cccccccccc} 
\hline
  $N_G$ & $L_2$ & $\mathcal{O}(L_2)$ & $L_2$ & $\mathcal{O}(L_2)$ &\\ 
\hline
  \multicolumn{3}{c}{Variable $v_y$} & \multicolumn{2}{c}{Variable $E_y$}  \\
\hline
  \multicolumn{5}{c}{$\mathcal{O}3$}    \\
\hline
100  & 8.2658E-02 &     &  1.1382E-02 &     \\ 
200  & 1.2933E-02 & 2.7 &  1.7733E-03 & 2.7 \\ 
400  & 1.6965E-03 & 2.9 &  2.3318E-04 & 2.9 \\ 
800  & 2.1272E-04 & 3.0 &  2.9843E-05 & 3.0 \\ 
\hline 
	\multicolumn{5}{c}{$\mathcal{O}4$}    \\
\hline
100  & 1.5487E-02 &     & 3.0653E-03 &     \\ 
200  & 5.1306E-04 & 4.9 & 9.7284E-05 & 5.0 \\ 
400  & 1.9320E-05 & 4.7 & 3.2709E-06 & 4.9 \\ 
800  & 9.3922E-07 & 4.4 & 1.8739E-07 & 4.1 \\ 
\hline 
  \multicolumn{5}{c}{$\mathcal{O}5$}    \\
\hline
100  & 5.9185E-03 &     & 7.2904E-04 &     \\ 
200  & 2.9331E-04 & 4.3 & 3.6769E-05 & 4.3 \\ 
300  & 4.2396E-05 & 4.8 & 5.3220E-06 & 4.8 \\ 
400  & 1.0396E-05 & 4.9 & 1.3026E-06 & 4.9 \\ 
\hline 
  \multicolumn{5}{c}{$\mathcal{O}6$}    \\
\hline
50   & 1.8839E-02 &     & 3.6056E-03 &     \\ 
100  & 6.2118E-04 & 4.9 & 1.1061E-04 & 5.0 \\ 
200  & 1.0376E-05 & 5.9 & 1.6986E-06 & 6.0 \\ 
300  & 8.4098E-07 & 6.2 & 1.5450E-07 & 5.9 \\ 
\hline 
\end{tabular} 
\end{center}
\label{tab.conv2}
\end{table} 

\subsubsection{Shock tube problems} 
\label{sec.shock.rrmhd}

\begin{figure}[!htbp]
\begin{center}
\begin{tabular}{lr}
\includegraphics[width=0.45\textwidth]{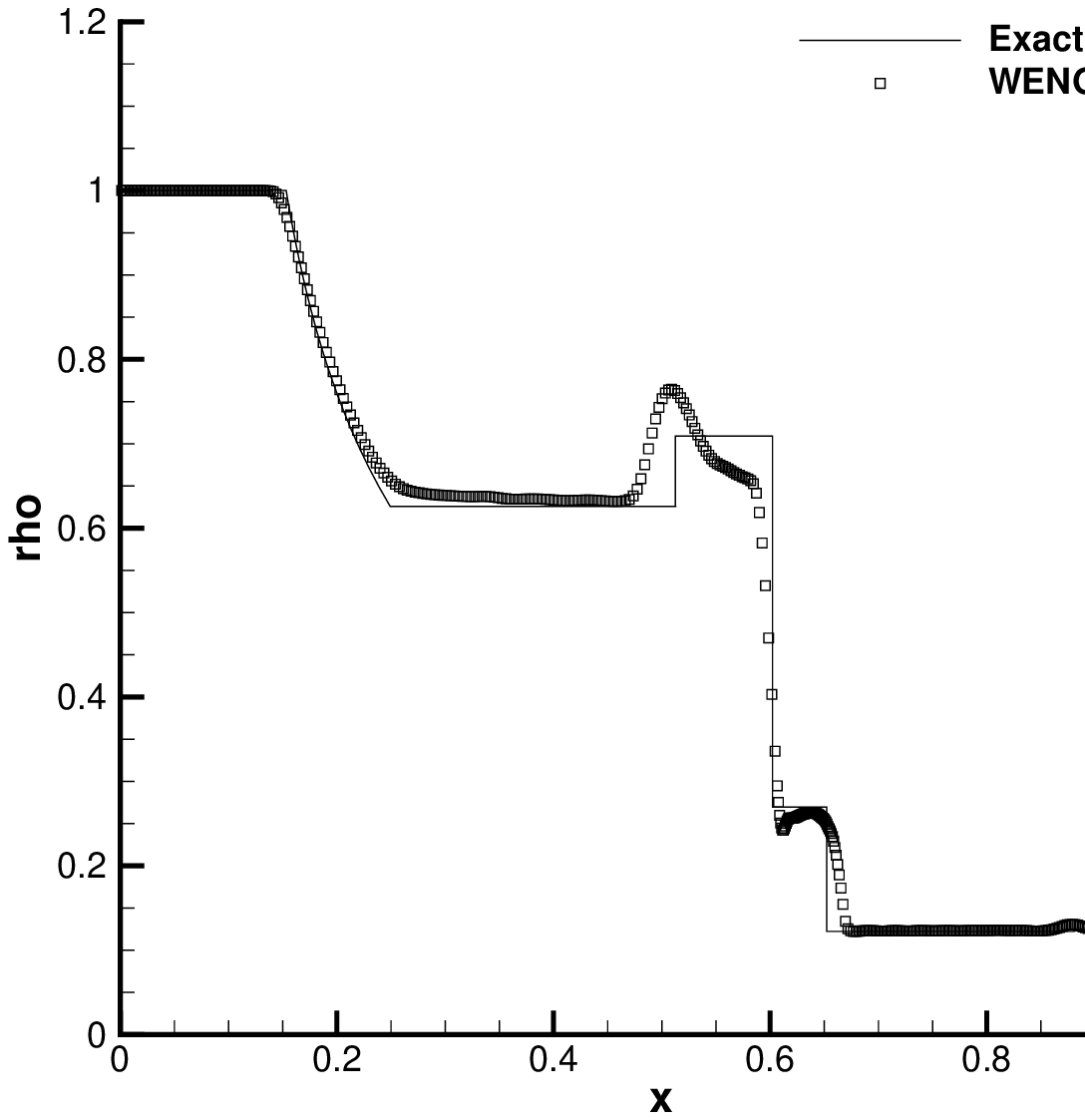} &
\includegraphics[width=0.45\textwidth]{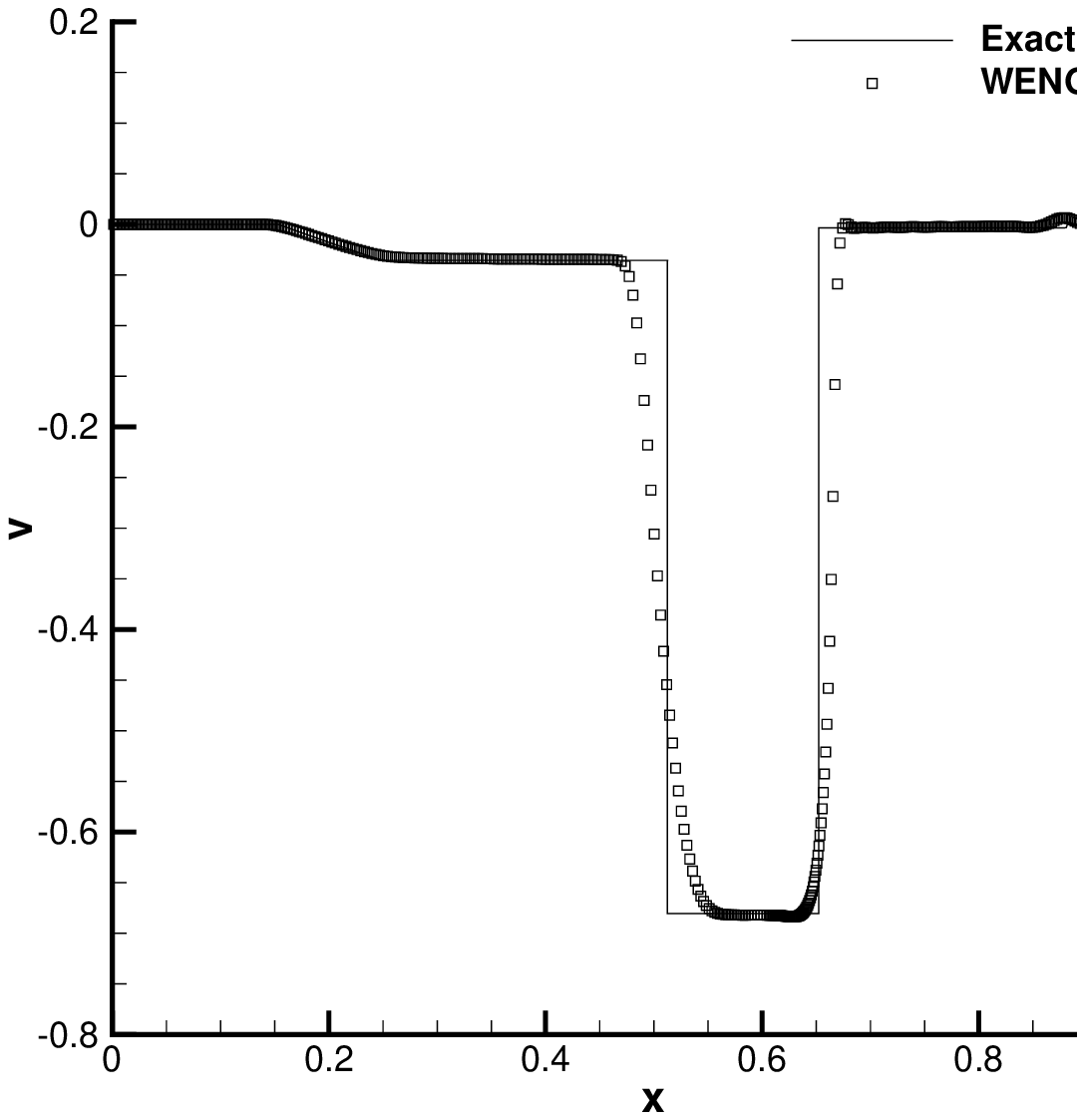} \\ 
\includegraphics[width=0.45\textwidth]{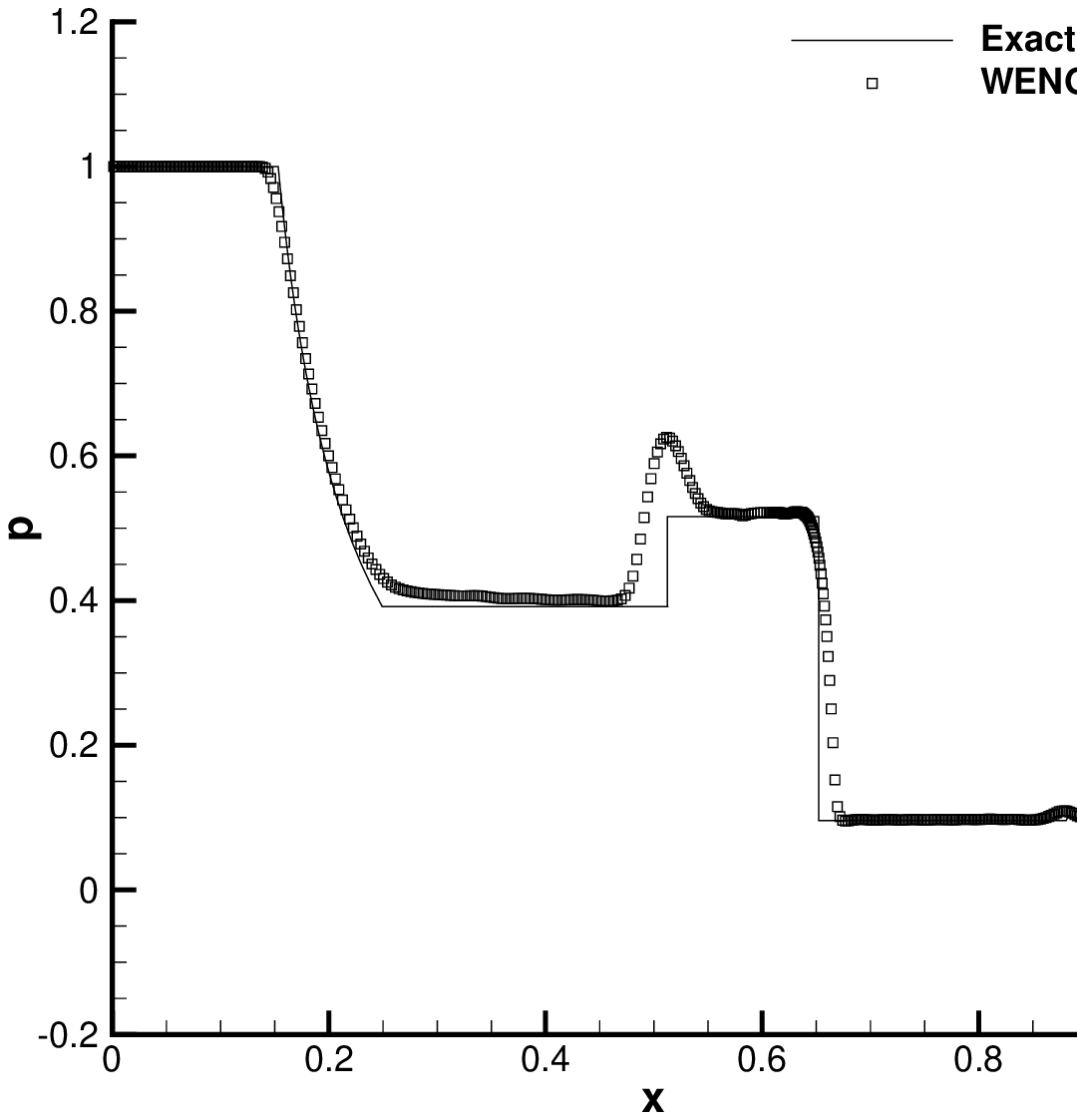} &
\includegraphics[width=0.45\textwidth]{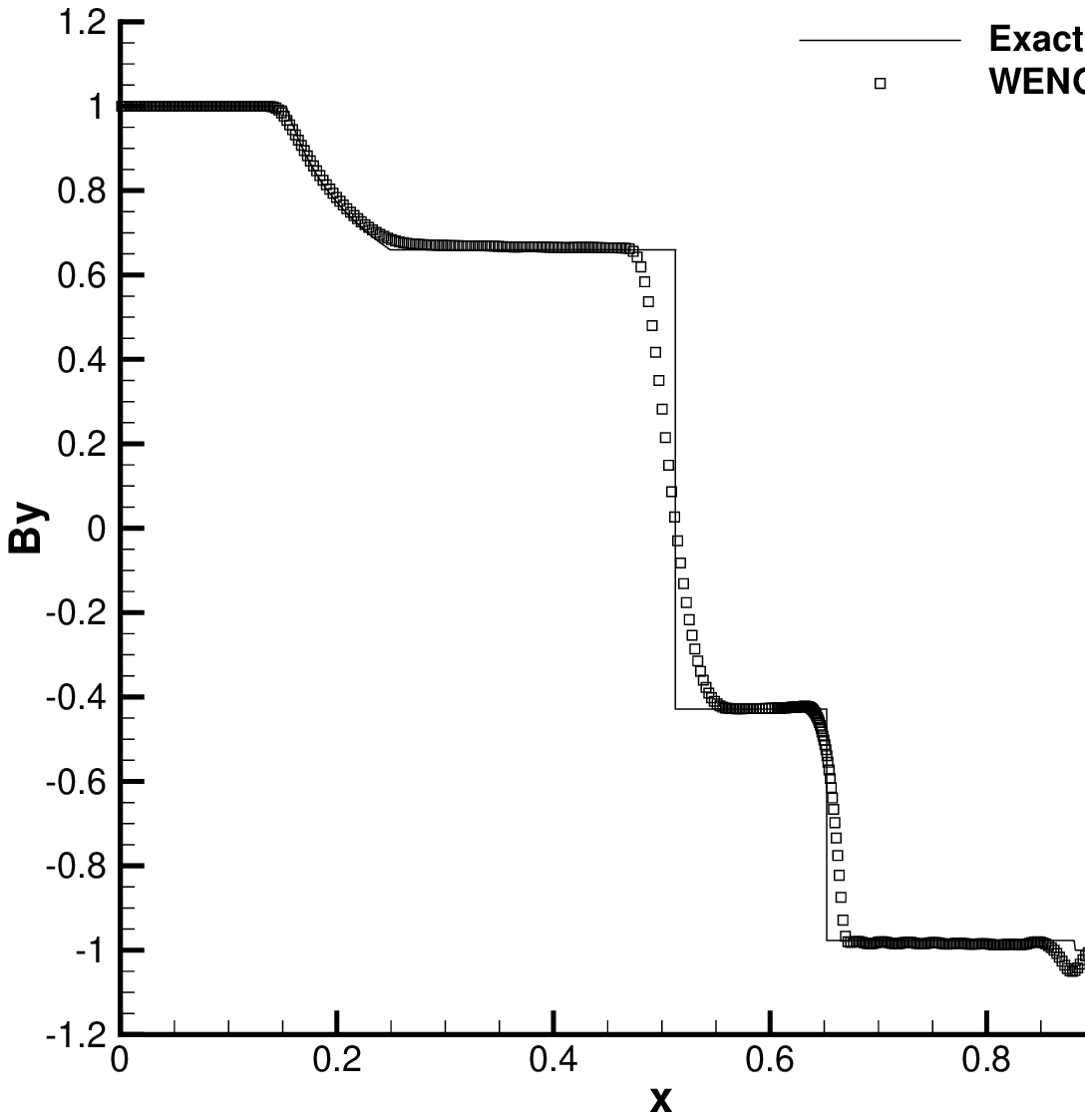} \\  
\end{tabular}
\caption{Numerical solution obtained for shock tube problem 1 ($\sigma=10^3$) of the resistive relativistic MHD equations using a third 
order Lagrangian one--step WENO finite volume scheme and exact solution of the ideal RMHD equations. Results are shown for density $\rho$, transverse velocity $v$, pressure $p$ and magnetic field $B_y$. }
\label{fig.rrmhd.rp1} 
\end{center}
\end{figure}

\begin{figure}[!htbp]
\begin{center}
\begin{tabular}{lr}
\includegraphics[width=0.45\textwidth]{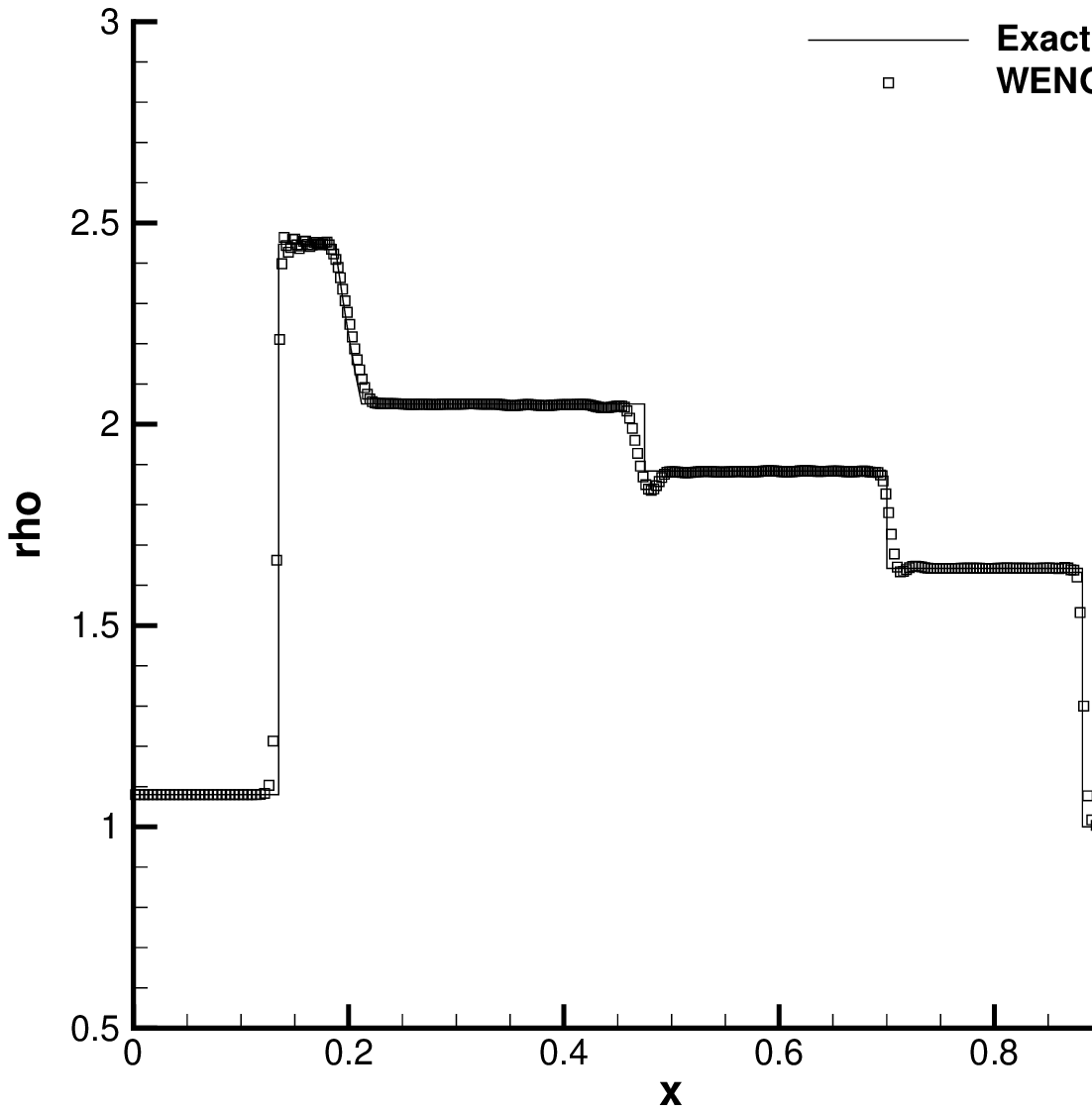} &
\includegraphics[width=0.45\textwidth]{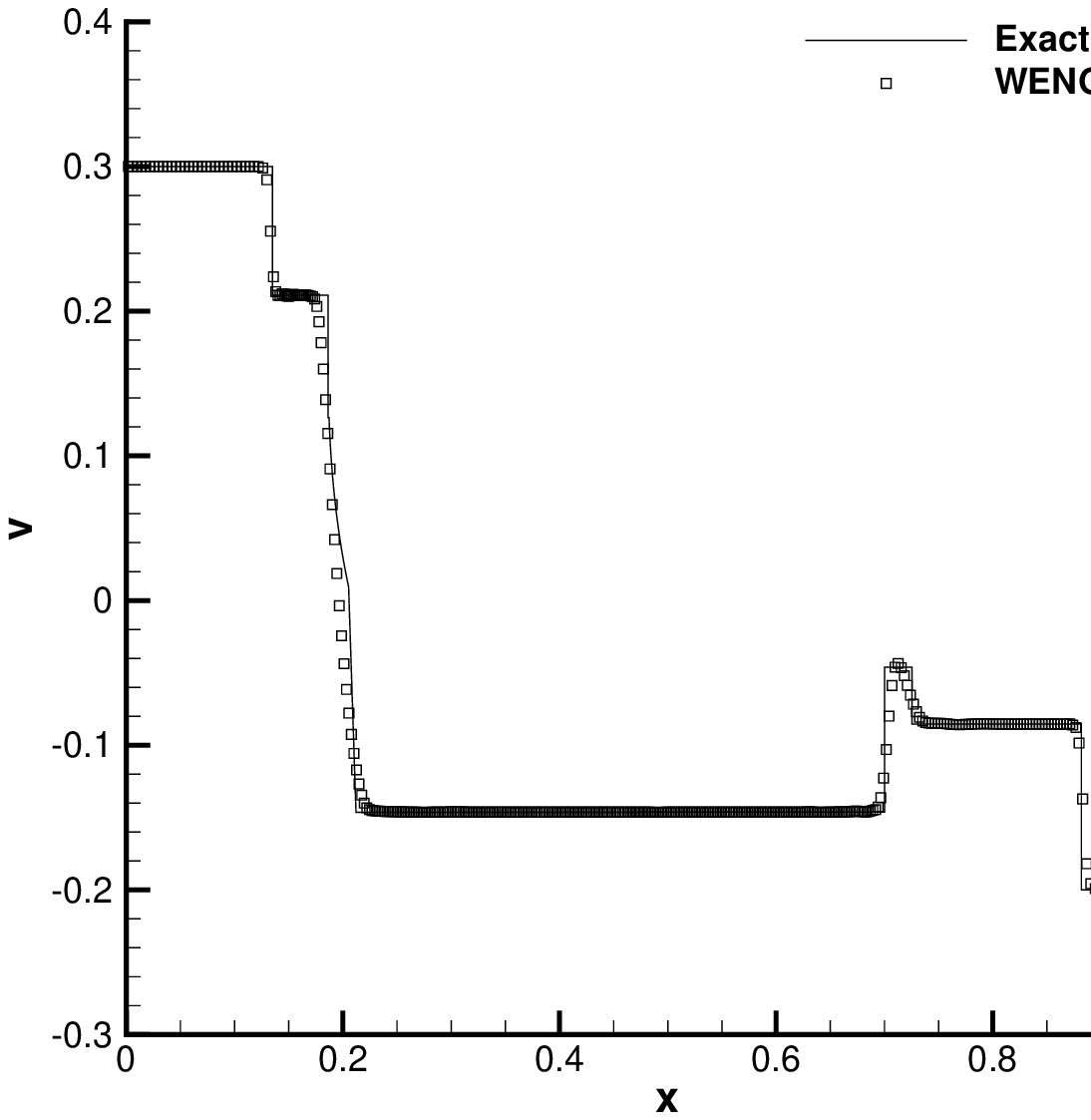} \\ 
\includegraphics[width=0.45\textwidth]{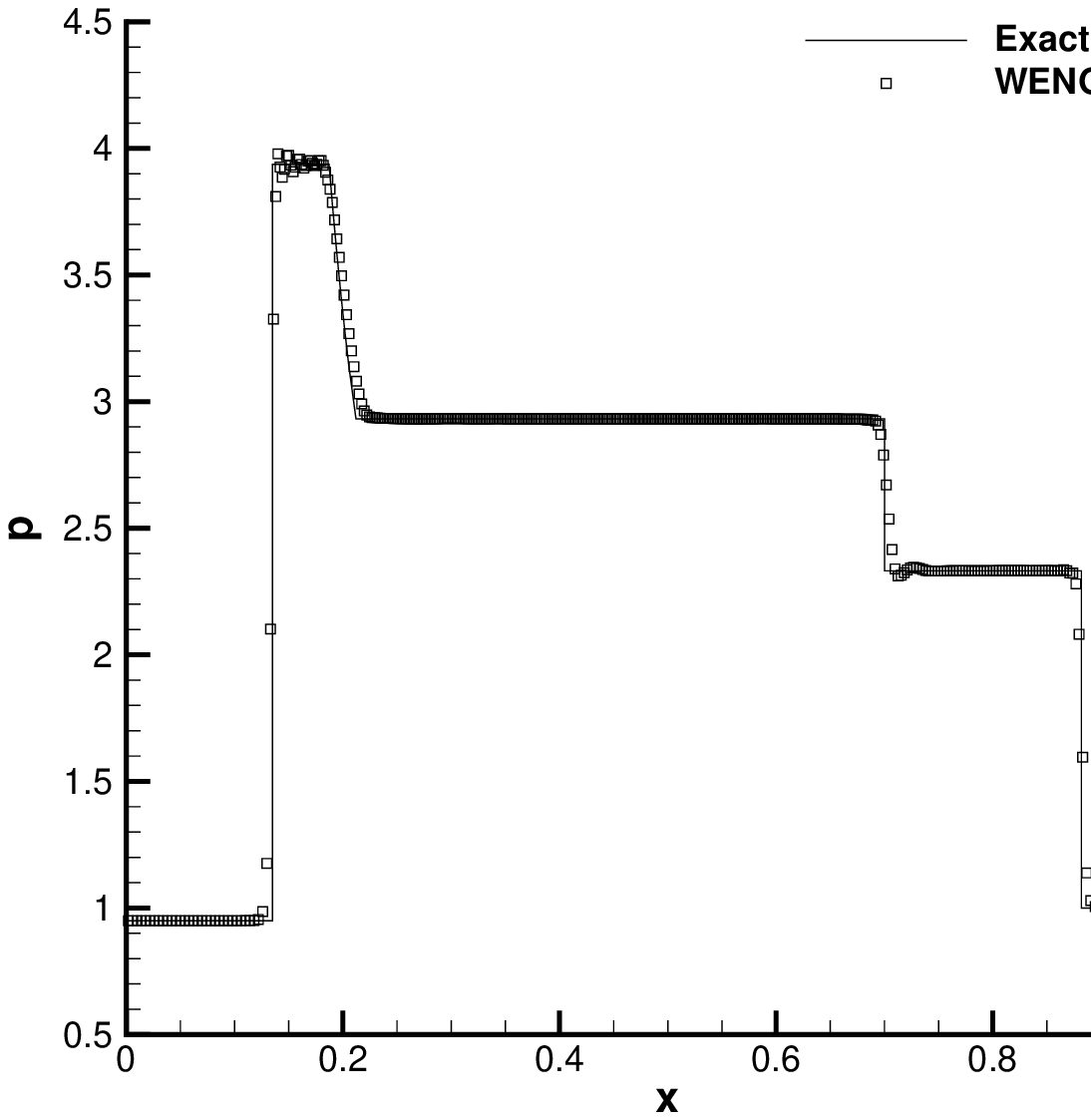} &
\includegraphics[width=0.45\textwidth]{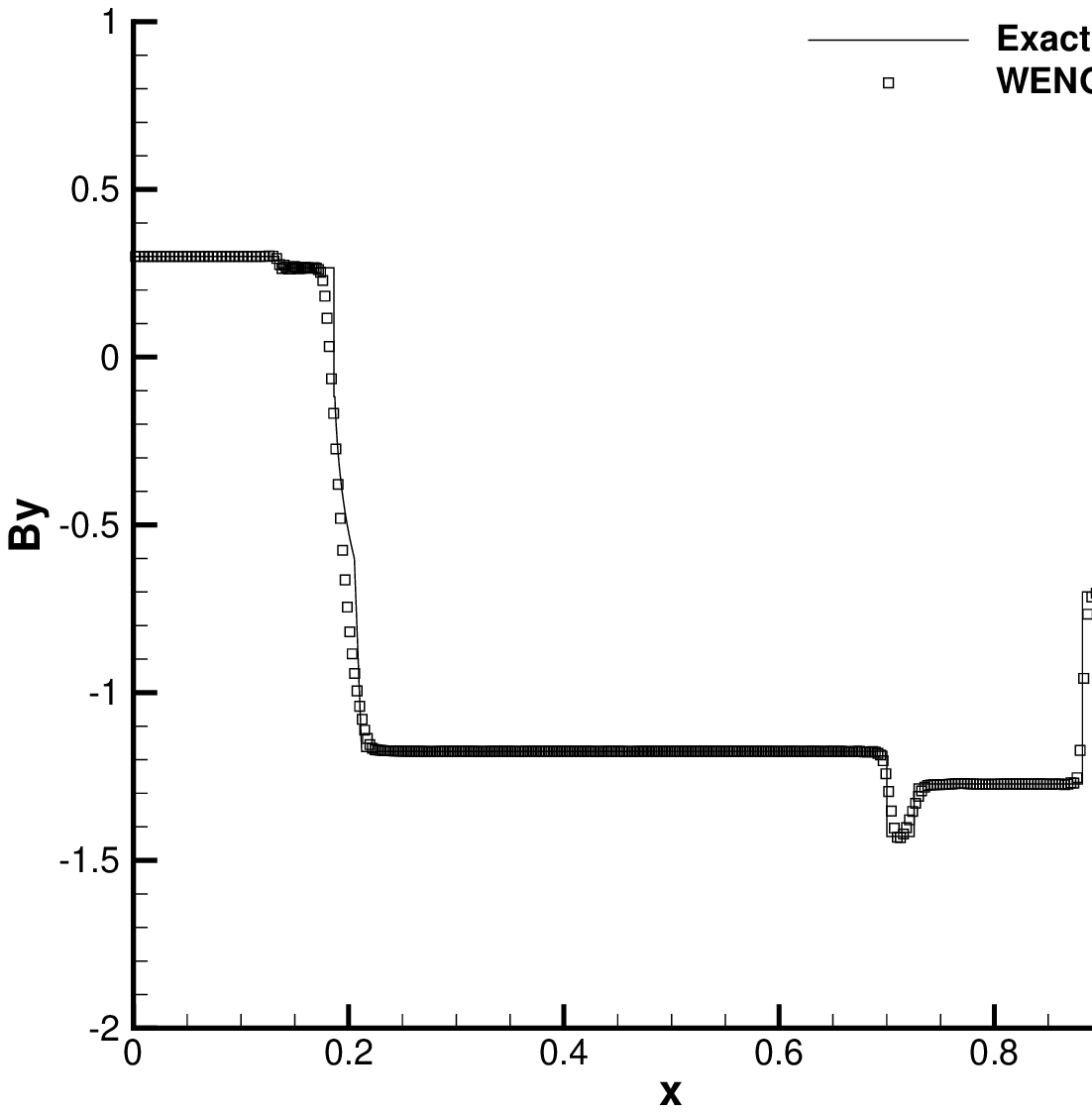} \\  
\end{tabular}
\caption{Numerical solution obtained for shock tube problem 2 ($\sigma=10^5$) of the resistive relativistic MHD equations using a third 
order Lagrangian one--step WENO finite volume scheme and exact solution of the ideal RMHD equations. Results are shown for density $\rho$, 
transverse velocity $v$, pressure $p$ and magnetic field $B_y$. }
\label{fig.rrmhd.rp2}
\end{center}
\end{figure}

In this section we solve two out of a series of test problems proposed by Balsara in \cite{BalsaraRMHD} for the ideal relativistic MHD equations. In particular, we solve the 
resistive RMHD equations with a large value for the conductivity $\sigma$ to validate the behaviour of our method in the presence of stiff source terms. The initial condition 
is given by two piecewise constant states separated by a discontinuity at $x=0$. The left and right values for the primitive variables are reported in Table \ref{tab.rmhd.ic}. 
Furthermore, we set $E^i = - \epsilon^{ijk} v^i B_k$, $\phi=\psi=q=0$ and $\gamma=\frac{5}{3}$. The conductivity in our test cases is chosen as $\sigma=10^3$ for the first
test problem and $\sigma=10^5$ for the second one. The computational domain is $\Omega = [-0.5;0.5]$ with Dirichlet boundaries consistent with the initial condition in $x$-direction. 
We use an initially equidistant grid with 400 points. The numerical results are shown together with the exact solution in Figures \ref{fig.rrmhd.rp1} and \ref{fig.rrmhd.rp2}. The exact solution is the one for the ideal RMHD equations, as published in \cite{GiacomazzoRezzolla}. The essential wave structures of the ideal RMHD Riemann problem can be noted. 

\begin{table}
\caption{Initial left (L) and right (R) states for the resistive relativistic MHD shock tube problems and final times $t_e$.}
\begin{center} 
\begin{tabular}{cccccccccc}
\hline
 Case  & $\rho$ & $p$ & $u$  & $v$ & $w$ & $B_y$ & $B_z$ & $B_x$ & $t_e$ \\
\hline 
1L & 1.0   & 1.0 & 0.0 & 0.0 & 0.0 &  1.0 & 0.0 & 0.5 & 0.4 \\
1R & 0.125 & 0.1 & 0.0 & 0.0 & 0.0 & -1.0 & 0.0 & 0.5 &     \\
\hline
2L & 1.08 & 0.95 & 0.4   & 0.3  & 0.2 & 0.3  & 0.3 & 2.0 & 0.55 \\
2R & 1.0  & 1.0  & -0.45 & -0.2 & 0.2 & -0.7 & 0.5 & 2.0 &      \\
\hline
\end{tabular}
\end{center}
\label{tab.rmhd.ic}
\end{table}

\subsection{Relativistic Radiation Hydrodynamics} 
\label{sec.radhydro}

As an additional test, we have considered the
solution of the special relativistic radiation
hydrodynamics equations. In the truncated
moment formalism introduced by \cite{Thorne1981}, it is
possible to write such equations, at least in the
optically thick regime, in the conservative 
form required by 
Eq.~\eqref{eqn.pde} (see \cite{Farris08}, \cite{Roedig2012} and
\cite{Zanotti2011}). 
In one spatial dimension, the vectors of the conservative 
variables in eqn. \eqref{eqn.pde}, of the fluxes and of the 
sources are respectively given by
\be
\Q=\left(\begin{array}{c}
D \\ S \\ \tau \\ S_{\rm r} \\ \tau_{\rm r}
\end{array}\right), \hspace{1cm}
\F = \left(\begin{array}{c}
v D \\
Z  \\
S \\
R_{\rm r} \\
S_{\rm r} \\
\end{array}\right) 
\hspace{1cm}
\S = \left(\begin{array}{c}
0 \\  
G_{\rm r} \\
 G^t_{\rm r} \\ 
 -G_{\rm r} \\ 
 -G^t_{\rm r}
\end{array}\right),
\label{eq:fluxes}
\ee
The first three equations  express
the usual conservation of mass, momentum and
energy of the fluid, and the corresponding 
conservative variables $(D,S,\tau)$ have the same definition as in
(\ref{eq:D})-(\ref{eq:U}), except for the fact that the
electromagnetic fields are zero.
The last two equations, on
the other hand, represent the time evolution of the flux
and of the energy density of the radiation field as measured
in the laboratory frame, with
\bea
S_{\rm r} & = & \frac{4}{3}E_{\rm r} W^2v +
W f_{\rm r}(1+ v^2) \,,
\label{eq:S_r} \\
\tau_{\rm r}& = &\frac{4}{3}E_{\rm r} W^2 + 2
W f_{\rm r} v - \frac{E_{\rm r}}{3} \,.
\label{eq:U_r}
\eea
The primitive variables of the radiation field are the
flux $f_{\rm r}$ and the energy density $E_{\rm r}$ as
measured in the comoving frame of the fluid, and they are
formally
related to the specific intensity of the radiation
$I_\nu$ \cite{Zanotti2011}. Finally, the quantities $Z$
and $R_{\rm r}$ entering (\ref{eq:fluxes}) are defined as 
\bea
Z &=& \omega W^2 v^2+ p \,, \\
R_{\rm r} &=& \frac{4}{3}E_{\rm r}W^2v^2 +2Wf_{\rm
  r}v+\frac{E_{\rm r}}{3}\,.
\eea
The sources of the radiation field $G^t_{\rm r}$ and
$G_{\rm r}$ depend on the physical interaction between
radiation and matter and can be written as
\bea
G^t_{\rm r} &=&\chi^t(E_{\rm r}-4\pi \tilde
B)W+(\chi^t+\chi^s)v f_{\rm r} \, \\
G_{\rm r} &=&\chi^t(E_{\rm r}-4\pi \tilde B)v W+(\chi^t+\chi^s)f_{\rm r} \,,
\eea
where $4\pi \tilde B=a_{\rm rad}T^4$ is the equilibrium black body
intensity, while $\chi^t$ and $\chi^s$ are the thermal and
the scattering opacity coefficients, respectively, which
are ultimately responsible for the stiffness of these
equations. The temperature $T$ of the fluid is computed
from the ideal-gas equation of state through the simple
relation $T=p/\rho$.
We also recall that, while the conversion
from the purely hydrodynamical conservative variables
to the corresponding primitive variables is not analytic
and it requires
the solution of an algebraic equation \cite{DelZanna2002}, the conversion
from the conservative radiation variables $(S_{\rm
  r},\tau_{\rm r})$ to the corresponding primitive
variables $(f_{\rm r},E_{\rm r})$ is just linear and
follows directly from (\ref{eq:S_r})-(\ref{eq:U_r}).

\begin{table}
\caption{Initial left (L) and right (R) states for the two
  relativistic radiation hydrodynamics 
shock tube problems considered.}
\begin{center} 
\begin{tabular}{llllllcc}
\hline
\vspace{0.05cm}
 Case  & $\rho$ & $p$ & $v$  & $E_{\rm r}$ &
 $\chi^t/\rho$ & $\gamma$ &$a_{\rm rad}$ \\
\hline 
\vspace{0.05cm}
1L & 1.0   & 4.0$\times10^{-3}$   & 0.2425   & 2.0$\times10^{-5}$ &  0.2 &5/3 &7.812$\times10^4$ \\
1R & 3.11   & 4.512$\times10^{-2}$ & 8.014$\times10^{-2}$ & 3.46$\times10^{-3}$ &  &   \\
\hline
\vspace{0.05cm}
2L & 1.0  & 60.0 & 0.995   & 2.0  & 0.3 & 2 & 1.543$\times10^{-7}$\\
2R & 8.0  & 2.34$\times10^3$  & 0.781 & 1.14$\times10^3$ &  &  &\\
\hline
\end{tabular}
\end{center}
\label{tab.rad.ic}
\end{table}

\begin{figure}[!htbp]
\begin{center}
\begin{tabular}{cc}
\includegraphics[width=0.45\textwidth]{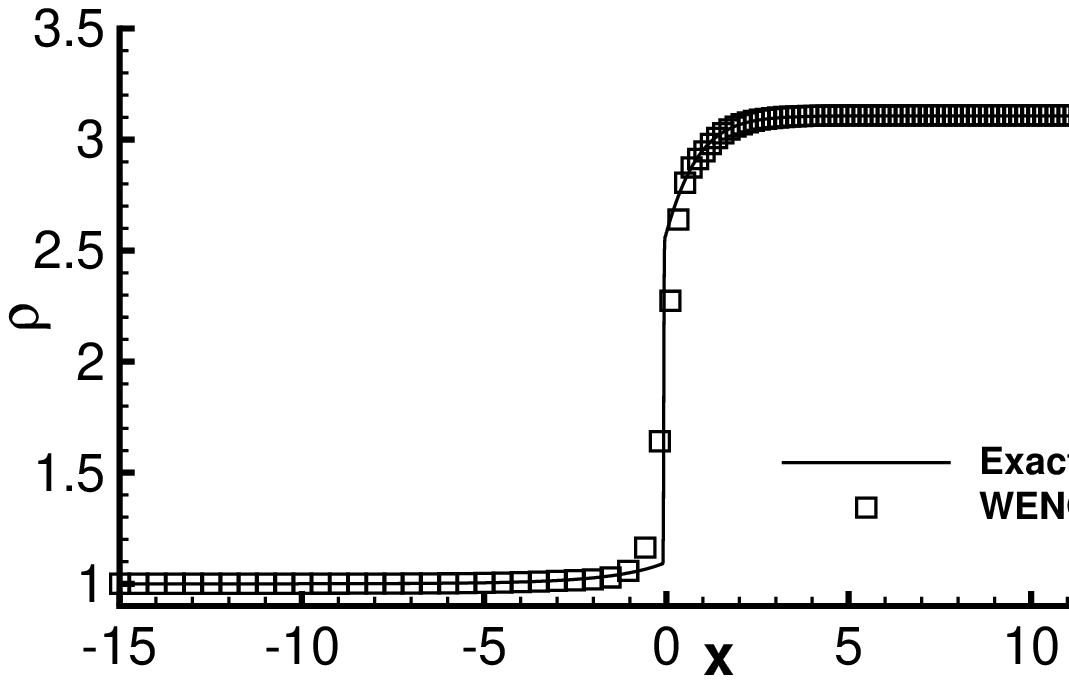}
& 
\includegraphics[width=0.45\textwidth]{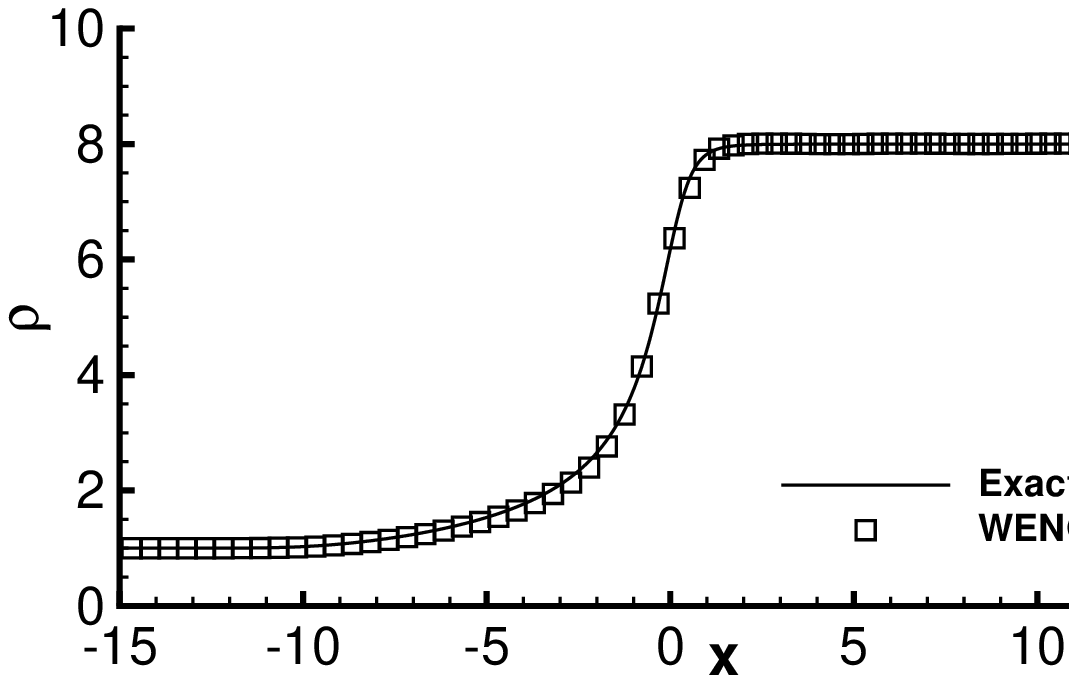} \\
\includegraphics[width=0.45\textwidth]{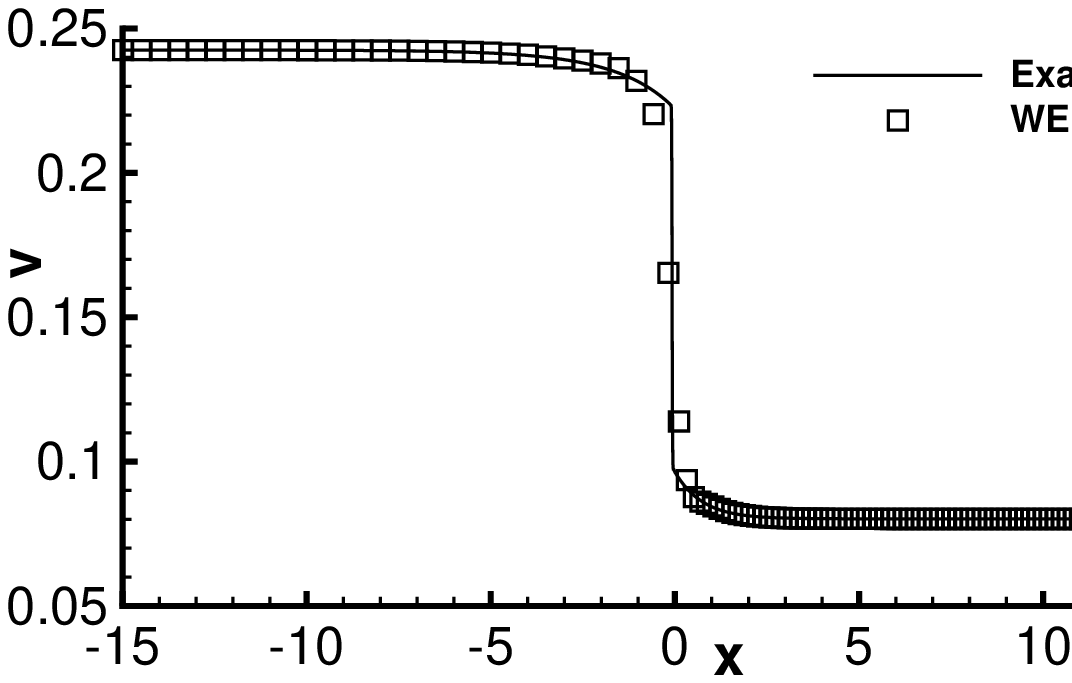}
& 
\includegraphics[width=0.45\textwidth]{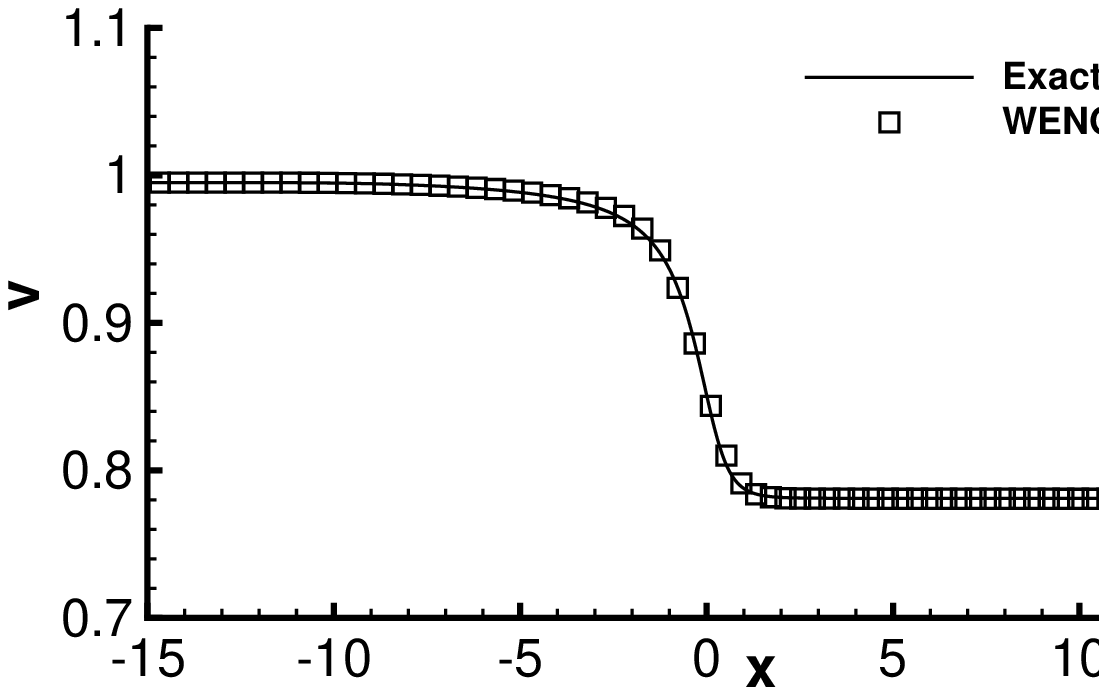} \\
\includegraphics[width=0.45\textwidth]{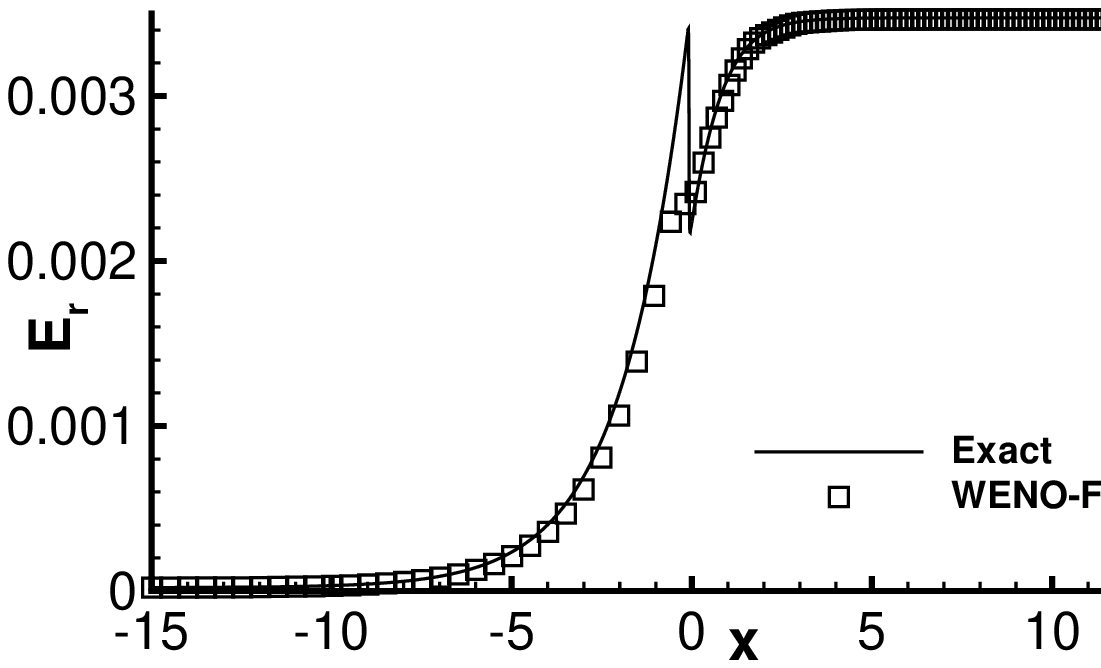}
& 
\includegraphics[width=0.45\textwidth]{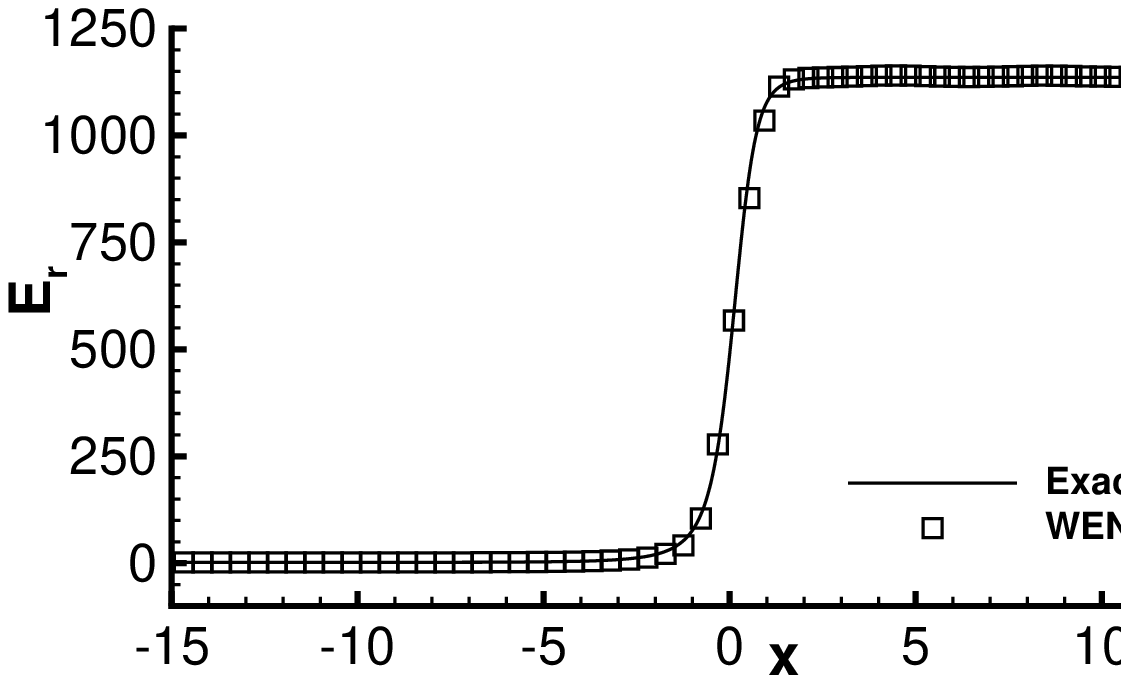}
\end{tabular}
\caption{
Numerical solution obtained for shock tube problem 1
(left panels) and for shock tube problem 2 (right panels)
of the relativistic radiation hydrodynamics equations
 using a third 
order Lagrangian one--step WENO finite volume scheme. 
Results are shown for density $\rho$, 
velocity $v$, and energy density of the radiation field
$E_{\rm r}$.
}
\label{fig.rad.shock.test2}
\end{center}
\end{figure}
In the verification of our numerical scheme, we
have considered two shock-tube tests, with initial
conditions reported in Tab.~(\ref{tab.rad.ic}).
The first test involves the propagation of a mildly
relativistic strong shock, while the second one generates
a smooth highly relativistic wave.
Each test is evolved in time until stationarity is reached. 
The semi-analytic solution that is used for comparison with the numerical
one has been obtained following the strategy described
by~\cite{Farris08}. The
scattering opacity coefficient $\chi^s$ has been set to
zero,
while the value of the thermal opacity coefficient $\chi^t$ is reported in
Table~\ref{tab.rad.ic} and it produces configurations
that are moderately
stiff. Figure~\ref{fig.rad.shock.test2} shows the 
comparison of the numerical solution with the exact one,
where we have adopted a third-order WENO reconstruction,
with $\textnormal{CFL}=0.4$ and 100 initially equidistant
grid points. 
\begin{figure}[!htbp]
\begin{center}
\begin{tabular}{lr}
\includegraphics[width=0.45\textwidth]{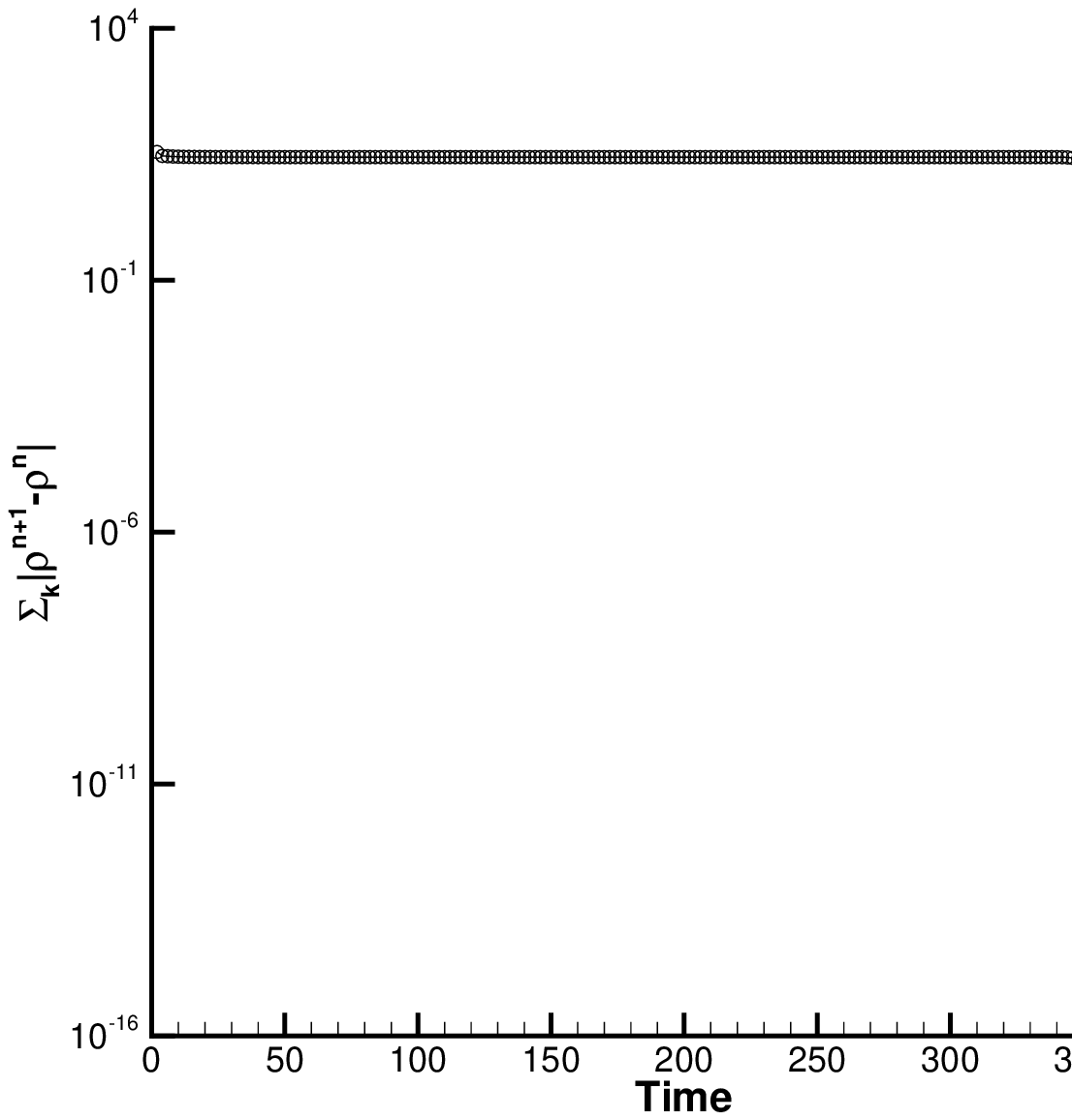} & 
\includegraphics[width=0.45\textwidth]{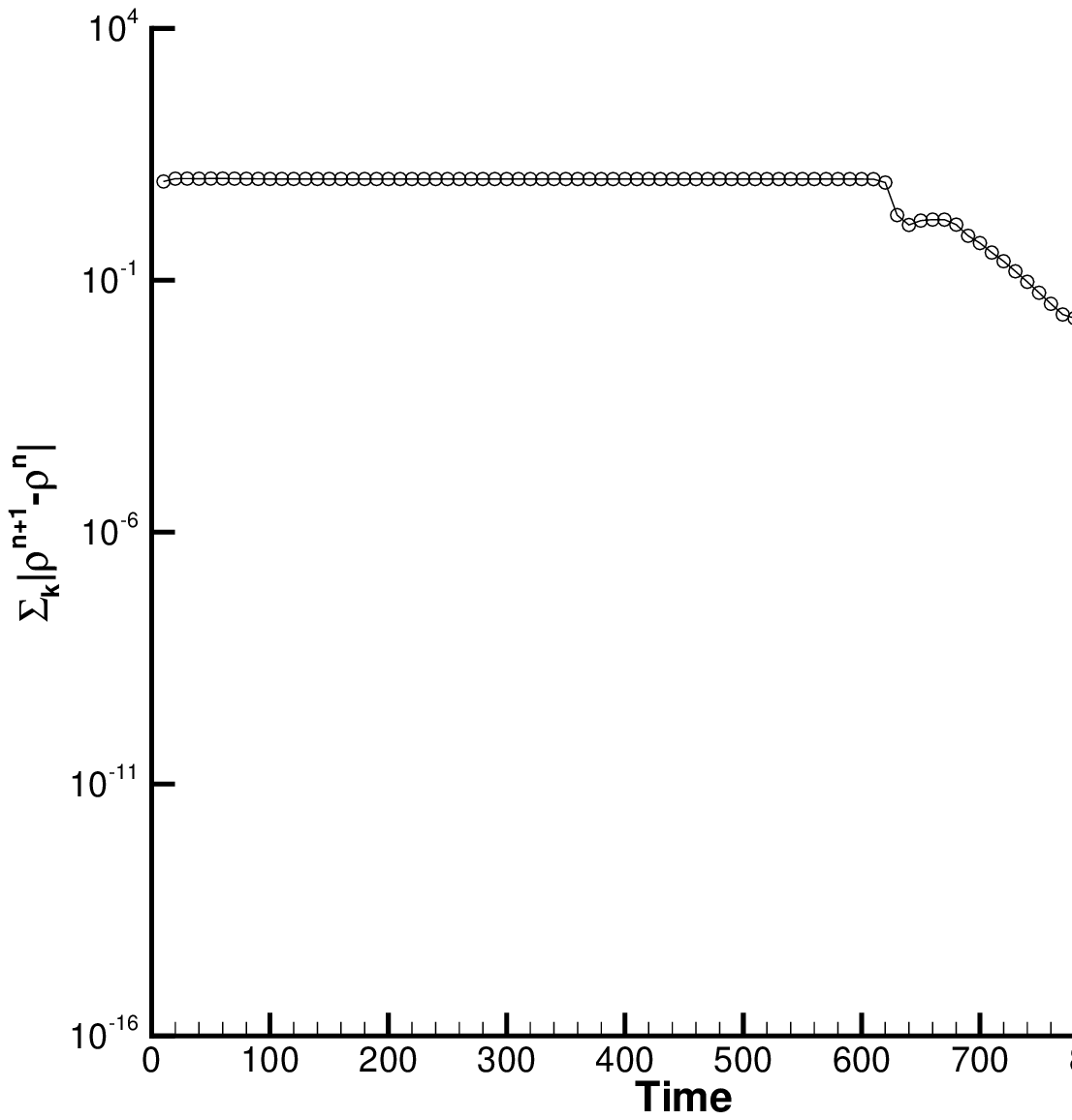} \\
\end{tabular}
\caption{
Evolution of the residual in the variable $\rho$ for test problem
1 (left) and test problem 2 (right).} 
\label{fig.residual}
\end{center}
\end{figure}
In Figure \ref{fig.residual} we show the evolution of the 
residual of the density $\rho$. It can be noted that initially the residual 
stagnates at a rather high level due to the presence of many transient
waves inside the computational domain. Once the transient waves have left, the 
residual drops very quickly to machine precision. These results confirm the ability 
of the new scheme in solving the relativistic radiation hydrodynamics equations.
It has furthermore been shown that the proposed high order numerical scheme is 
able to simulate correctly, both, time dependent problems as well as steady 
state problems. 
\section{Conclusions} 
\label{sec.conclusions} 

We have developed a new high order one-step Arbitrary-Lagrangian-Eulerian (ALE) WENO finite volume scheme for the solution of nonlinear 
systems of hyperbolic balance laws with stiff source terms. The presented approach has been validated against exact reference solutions 
available for smooth and discontinuous solutions of three different hyperbolic systems, namely the Euler equations of gas dynamics, the 
resistive relativistic MHD equations and the relativistic radiation-hydrodynamics equations. In all cases the algorithm was
found to be very robust and at the same time very accurate. To our knowledge, it is the first time that such high orders 
have been reached with ALE-type finite volume schemes. In the near future we plan to extend the schemes presented in this 
article to structured and unstructured meshes in multiple space dimensions. Again, the building blocks will be a high order 
WENO reconstruction \cite{DumbserKaeser06b,DumbserKaeser07} and a local space-time DG predictor 
\cite{DumbserEnauxToro,DumbserZanotti,HidalgoDumbser}.  

%%%% Acknowledgments %%%%%%%%
\section*{Acknowledgments}
The authors would like to thank B. Giacomazzo and L. Rezzolla for providing the exact Riemann solvers for ideal RMHD. 
The presented research has been financed by the European Research Council under the European Union's Seventh Framework 
Programme (FP7/2007-2013) under the research project \textit{STiMulUs}, ERC Grant agreement no. 278267. 
The stay of A.U. at the University of Trento was co--financed by the ERASMUS mundus external cooperation 
window \textit{Bridging the Gap - Mobility lot 12} and was made possible thanks to Prof. Sarantuya Tsedendamba from the 
Mongolian University of Science and Technology.  

The authors would like to thank the two anonymous referees for their constructive suggestions and comments that greatly 
helped to improve the quality and the clarity of this article.  

\bibliography{Lagrange1D}

\begin{thebibliography}{10}

\bibitem{BalsaraRMHD}
D.~Balsara.
\newblock Total variation diminishing scheme for relativistic
  magnetohydrodynamics.
\newblock {\em The Astrophysical Journal Supplement Series}, 132:83--101, 2001.

\bibitem{balsarashu}
D.~Balsara and {C.W.} Shu.
\newblock Monotonicity preserving weighted essentially non-oscillatory schemes
  with increasingly high order of accuracy.
\newblock {\em Journal of Computational Physics}, 160:405--452, 2000.

\bibitem{Artzi}
M.~Ben-Artzi and J.~Falcovitz.
\newblock A second-order godunov-type scheme for compressible fluid dynamics.
\newblock {\em Journal of Computational Physics}, 55:1--32, 1984.

\bibitem{Benson1992}
D.J. Benson.
\newblock Computational methods in lagrangian and eulerian hydrocodes.
\newblock {\em Computer Methods in Applied Mechanics and Engineering},
  99:235--394, 1992.

\bibitem{Raviart.GRP.2}
A.~Bourgeade, P.~LeFloch, and P.A. Raviart.
\newblock An asymptotic expansion for the solution of the generalized {Riemann}
  problem. {Part II}: application to the gas dynamics equations.
\newblock {\em Annales de l'institut Henri Poincar\'e (C) Analyse non
  lin\'eaire}, 6:437--480, 1989.

\bibitem{Caramana1998}
E.J. Caramana, D.E. Burton, M.J. Shashkov, and P.P. Whalen.
\newblock The construction of compatible hydrodynamics algorithms utilizing
  conservation of total energy.
\newblock {\em Journal of Computational Physics}, 146:227--262, 1998.

\bibitem{Carre2009}
G.~Carr\'e, S.~Del Pino, B.~Despr\'es, and E.~Labourasse.
\newblock A cell-centered lagrangian hydrodynamics scheme on general
  unstructured meshes in arbitrary dimension.
\newblock {\em Journal of Computational Physics}, 228:5160--5183, 2009.

\bibitem{Casulli1990}
V.~Casulli.
\newblock Semi-implicit finite difference methods for the two-dimensional
  shallow water equations.
\newblock {\em Journal of Computational Physics}, 86:56--74, 1990.

\bibitem{CasulliCheng1992}
V.~Casulli and R.T. Cheng.
\newblock Semi-implicit finite difference methods for three-dimensional shallow
  water flow.
\newblock {\em International Journal of Numerical Methods in Fluids},
  15:629--648, 1992.

\bibitem{chengshu1}
J.~Cheng and C.W. Shu.
\newblock {A high order ENO conservative Lagrangian type scheme for the
  compressible Euler equations}.
\newblock {\em Journal of Computational Physics}, 227:1567--1596, 2007.

\bibitem{chengshu3}
J.~Cheng and C.W. Shu.
\newblock {A cell-centered Lagrangian scheme with the preservation of symmetry
  and conservation properties for compressible fluid flows in two-dimensional
  cylindrical geometry}.
\newblock {\em Journal of Computational Physics}, 229:7191--7206, 2010.

\bibitem{chengshu4}
J.~Cheng and C.W. Shu.
\newblock {Improvement on spherical symmetry in two-dimensional cylindrical
  coordinates for a class of control volume Lagrangian schemes}.
\newblock {\em Communications in Computational Physics}, 11:1144--1168, 2012.

\bibitem{CIR}
R.~Courant, E.~Isaacson, and M.~Rees.
\newblock On the solution of nonlinear hyperbolic differential equations by
  finite differences.
\newblock {\em Comm. Pure Appl. Math.}, 5:243--255, 1952.

\bibitem{Dedneretal}
A.~Dedner, F.~Kemm, D.~Kr\"oner, C.-D. Munz, T.~Schnitzer, and M.~Wesenberg.
\newblock Hyperbolic divergence cleaning for the {MHD} equations.
\newblock {\em Journal of Computational Physics}, 175:645--673, 2002.

\bibitem{DelZanna2002}
L.~{Del Zanna} and N.~{Bucciantini}.
\newblock {An efficient shock-capturing central-type scheme for
  multidimensional relativistic flows. I. Hydrodynamics}.
\newblock {\em Astron. Astroph.}, 390:1177--1186, August 2002.

\bibitem{delZanna2007}
L.~{Del Zanna}, O.~{Zanotti}, N.~{Bucciantini}, and P.~{Londrillo}.
\newblock {ECHO: a Eulerian conservative high-order scheme for general
  relativistic magnetohydrodynamics and magnetodynamics}.
\newblock {\em Astronomy \& Astrophysics}, 473:11--30, October 2007.

\bibitem{Despres2005}
B.~Despr\'es and C.~Mazeran.
\newblock Lagrangian gas dynamics in two-dimensions and lagrangian systems.
\newblock {\em Archive for Rational Mechanics and Analysis}, 178:327--372,
  2005.

\bibitem{ADERNSE}
M.~Dumbser.
\newblock Arbitrary high order {PNPM} schemes on unstructured meshes for the
  compressible {Navier--Stokes} equations.
\newblock {\em Computers \& Fluids}, 39:60--76, 2010.

\bibitem{Dumbser2008}
M.~Dumbser, D.S. Balsara, E.F. Toro, and C.D. Munz.
\newblock A unified framework for the construction of one-step finite-volume
  and discontinuous {Galerkin} schemes.
\newblock {\em Journal of Computational Physics}, 227:8209--–8253, 2008.

\bibitem{DumbserEnauxToro}
M.~Dumbser, C.~Enaux, and E.F. Toro.
\newblock Finite volume schemes of very high order of accuracy for stiff
  hyperbolic balance laws.
\newblock {\em Journal of Computational Physics}, 227:3971--–4001, 2008.

\bibitem{DumbserKaeser06b}
M.~Dumbser and M.~K\"aser.
\newblock Arbitrary high order non-oscillatory finite volume schemes on
  unstructured meshes for linear hyperbolic systems.
\newblock {\em Journal of Computational Physics}, 221:693--723, 2007.

\bibitem{DumbserKaeser07}
M.~Dumbser, M.~K\"aser, V.A Titarev, and E.F. Toro.
\newblock Quadrature-free non-oscillatory finite volume schemes on unstructured
  meshes for nonlinear hyperbolic systems.
\newblock {\em Journal of Computational Physics}, 226:204--243, 2007.

\bibitem{OsherUniversal}
M.~Dumbser and E.F. Toro.
\newblock On universal {Osher}--type schemes for general nonlinear hyperbolic
  conservation laws.
\newblock {\em Communications in Computational Physics}, 10:635--671, 2011.

\bibitem{DumbserZanotti}
M.~Dumbser and O.~Zanotti.
\newblock Very high order {PNPM} schemes on unstructured meshes for the
  resistive relativistic {MHD} equations.
\newblock {\em Journal of Computational Physics}, 228:6991--7006, 2009.

\bibitem{Farris08}
B.~D. {Farris}, T.~K. {Li}, Y.~T. {Liu}, and S.~L. {Shapiro}.
\newblock {Relativistic radiation magnetohydrodynamics in dynamical spacetimes:
  Numerical methods and tests}.
\newblock {\em Phys. Rev. D}, 78(2):024023, July 2008.

\bibitem{SPHWeirFlow}
A.~Ferrari.
\newblock {SPH simulation of free surface flow over a sharp-crested weir}.
\newblock {\em Advances in Water Resources}, 33:270--276, 2010.

\bibitem{SPHLagrange}
A.~Ferrari, M.~Dumbser, E.F. Toro, and A.~Armanini.
\newblock {A New Stable Version of the SPH Method in Lagrangian Coordinates}.
\newblock {\em Communications in Computational Physics}, 4:378--404, 2008.

\bibitem{SPH3D}
A.~Ferrari, M.~Dumbser, E.F. Toro, and A.~Armanini.
\newblock {A new 3D parallel SPH scheme for free surface flows}.
\newblock {\em Computers \& Fluids}, 38:1203--1217, 2009.

\bibitem{Dambreak3D}
A.~Ferrari, L.~Fraccarollo, M.~Dumbser, E.F. Toro, and A.~Armanini.
\newblock Three--dimensional flow evolution after a dambreak.
\newblock {\em Journal of Fluid Mechanics}, 663:456--477, 2010.

\bibitem{Raviart.GRP.1}
P.~Le Floch and P.A. Raviart.
\newblock An asymptotic expansion for the solution of the generalized {Riemann}
  problem. {Part I}: General theory.
\newblock {\em Annales de l'institut Henri Poincar\'e (C) Analyse non
  lin\'eaire}, 5:179--207, 1988.

\bibitem{GiacomazzoRezzolla}
B.~Giacomazzo and L.~Rezzolla.
\newblock The exact solution of the {Riemann} problem in relativistic
  magnetohydrodynamics.
\newblock {\em Journal of Fluid Mechanics}, 562:223--259, 2006.

\bibitem{eno}
A.~Harten, B.~Engquist, S.~Osher, and S.~Chakravarthy.
\newblock Uniformly high order essentially non-oscillatory schemes, {III}.
\newblock {\em Journal of Computational Physics}, 71:231--303, 1987.

\bibitem{HidalgoDumbser}
A.~Hidalgo and M.~Dumbser.
\newblock {ADER} schemes for nonlinear systems of stiff
  advection–diffusion–reaction equations.
\newblock {\em Journal of Scientific Computing}, 48:173--189, 2011.

\bibitem{Hirt1974}
C.~Hirt, A.~Amsden, and J.~Cook.
\newblock An arbitrary lagrangian–eulerian computing method for all flow
  speeds.
\newblock {\em Journal of Computational Physics}, 14:227–253, 1974.

\bibitem{HuiCoord}
W.H. Hui.
\newblock {The unified coordinate system in computational fluid dynamics}.
\newblock {\em Communications in Computational Physics}, 2:577--610, 2007.

\bibitem{Jua2011}
Zupeng Jia and Shudao Zhang.
\newblock A new high-order discontinuous galerkin spectral finite element
  method for lagrangian gas dynamics in two-dimensions.
\newblock {\em Journal of Computational Physics}, 230:2496--2522, 2011.

\bibitem{shu_efficient_weno}
{G.-S.} Jiang and {C.W.} Shu.
\newblock Efficient implementation of weighted {ENO} schemes.
\newblock {\em Journal of Computational Physics}, pages 202--228, 1996.

\bibitem{Komissarov2007}
S.~S. {Komissarov}.
\newblock {Multidimensional numerical scheme for resistive relativistic
  magnetohydrodynamics}.
\newblock {\em Mon. Not. Roy. Astr. Soc.}, 382:995--1004, December 2007.

\bibitem{LentineEtAl2011}
M.~Lentine, J\'on~T\'omas Gr\'etarsson, and R.~Fedkiw.
\newblock An unconditionally stable fully conservative semi-lagrangian method.
\newblock {\em Journal of Computational Physics}, 230:2857--2879, 2011.

\bibitem{chengshu2}
W.~Liu, J.~Cheng, and C.W. Shu.
\newblock {High order conservative Lagrangian schemes with Lax–Wendroff type
  time discretization for the compressible Euler equations}.
\newblock {\em Journal of Computational Physics}, 228:8872--8891, 2009.

\bibitem{Maire2011}
P.-H. Maire.
\newblock A high-order one-step sub-cell force-based discretization for
  cell-centered lagrangian hydrodynamics on polygonal grids.
\newblock {\em Computers and Fluids}, 46(1):341--347, 2011.

\bibitem{Maire2010}
P.-H. Maire.
\newblock A unified sub-cell force-based discretization for cell-centered
  lagrangian hydrodynamics on polygonal grids.
\newblock {\em International Journal for Numerical Methods in Fluids},
  65:1281–1294, 2011.

\bibitem{Maire2007}
P.H. Maire, R.~Abgrall, J.~Breil, and J.~Ovadia.
\newblock A cell-centered lagrangian scheme for two-dimensional compressible
  flow problems.
\newblock {\em SIAM Journal on Scientific Computing}, 29:1781–1824, 2007.

\bibitem{Dimitri2011}
D.~J. Mavriplis and C.~R. Nastase.
\newblock On the geometric conservation law for high order discontinuous
  galerkin discretizations on dynamically deforming meshes.
\newblock {\em Journal of Computational Physics}, 230:4285--4300, 2011.

\bibitem{Monaghan1994}
J.J. Monaghan.
\newblock Simulating free surface flows with {SPH}.
\newblock {\em Journal of Computational Physics}, 110:399--406, 1994.

\bibitem{munz94}
C.D. Munz.
\newblock {On Godunov--type schemes for Lagrangian gas dynamics}.
\newblock {\em SIAM Journal on Numerical Analysis}, 31:17--42, 1994.

\bibitem{Palenzuela2009}
C.~Palenzuela, L.~Lehner, O.~Reula, and L.~Rezzolla.
\newblock Beyond ideal {MHD}: towards a more realistic modeling of relativistic
  astrophysical plasmas.
\newblock {\em Mon. Not. R. Astron. Soc.}, 2009.

\bibitem{Peery2000}
J.S. Peery and D.E. Carroll.
\newblock Multi-material ale methods in unstructured grids,.
\newblock {\em Computer Methods in Applied Mechanics and Engineering},
  187:591--619, 2000.

\bibitem{QuiShu2011}
Jing-Mei Qiu and Chi-Wang Shu.
\newblock Conservative high order semi-lagrangian finite difference weno
  methods for advection in incompressible flow.
\newblock {\em Journal of Computational Physics}, 230:863--889, 2011.

\bibitem{roe}
{P.L.} Roe.
\newblock Approximate {Riemann} solvers, parameter vectors, and difference
  schemes.
\newblock {\em Journal of Computational Physics}, 43:357--372, 1981.

\bibitem{Roedig2012}
C.~{Roedig}, O.~{Zanotti}, and D.~{Alic}.
\newblock {General relativistic radiation hydrodynamics of accretion flows: II.
  Treating stiff source terms and exploring physical limitations}.
\newblock {\em ArXiv e-prints: astro-ph.HE/1206.6662}, June 2012.

\bibitem{Rusanov:1961a}
V.~V. Rusanov.
\newblock {Calculation of Interaction of Non--Steady Shock Waves with
  Obstacles}.
\newblock {\em J. Comput. Math. Phys. USSR}, 1:267--279, 1961.

\bibitem{Smith1999}
R.W.Smith.
\newblock {AUSM(ALE)}: a geometrically conservative arbitrary
  lagrangian–eulerian flux splitting scheme.
\newblock {\em Journal of Computational Physics}, 150:268–286, 1999.

\bibitem{stroud}
{A.H.} Stroud.
\newblock {\em Approximate Calculation of Multiple Integrals}.
\newblock Prentice-Hall Inc., Englewood Cliffs, New Jersey, 1971.

\bibitem{Thorne1981}
K.~S. {Thorne}.
\newblock {Relativistic radiative transfer - Moment formalisms}.
\newblock {\em Mon. Not. R. Astron. Soc.}, 194:439--473, February 1981.

\bibitem{toro3}
{V.A.} Titarev and {E.F.} Toro.
\newblock {ADER}: Arbitrary high order {Godunov} approach.
\newblock {\em Journal of Scientific Computing}, 17:609--618, 2002.

\bibitem{titarevtoro}
V.A. Titarev and E.F. Toro.
\newblock {ADER} schemes for three-dimensional nonlinear hyperbolic systems.
\newblock {\em Journal of Computational Physics}, 204:715--736, 2005.

\bibitem{Toro:2006a}
E.~F. Toro and V.~A. Titarev.
\newblock {Derivative Riemann solvers for systems of conservation laws and ADER
  methods}.
\newblock {\em Journal of Computational Physics}, 212(1):150--165, 2006.

\bibitem{toro-book}
E.F. Toro.
\newblock {\em {Riemann} Solvers and Numerical Methods for Fluid Dynamics}.
\newblock Springer, third edition, 2009.

\bibitem{toro4}
{E.F.} Toro and {V. A.} Titarev.
\newblock Solution of the generalized {Riemann} problem for advection-reaction
  equations.
\newblock {\em Proc. Roy. Soc. London}, pages 271--281, 2002.

\bibitem{titarevtoro2}
E.F. Toro and V.A. Titarev.
\newblock {ADER} schemes for scalar hyperbolic conservation laws with source
  terms in three space dimensions.
\newblock {\em Journal of Computational Physics}, 202:196--215, 2005.

\bibitem{spacetimedg1}
J.~J.~W. van~der Vegt and H.~van~der Ven.
\newblock Space–-time discontinuous {Galerkin} finite element method with
  dynamic grid motion for inviscid compressible flows {I}. general formulation.
\newblock {\em Journal of Computational Physics}, 182:546–--585, 2002.

\bibitem{spacetimedg2}
H.~van~der Ven and J.~J.~W. van~der Vegt.
\newblock Space-–time discontinuous {Galerkin} finite element method with
  dynamic grid motion for inviscid compressible flows {II}. efficient flux
  quadrature.
\newblock {\em Comput. Methods Appl. Mech. Engrg.}, 191:4747–--4780, 2002.

\bibitem{Neumann1950}
J.~von Neumann and R.D. Richtmyer.
\newblock A method for the calculation of hydrodynamics shocks.
\newblock {\em Journal of Applied Physics}, 21:232--237, 1950.

\bibitem{RMHD}
L.~Del Zanna, N.~Bucciantini, and P.~Londrillo.
\newblock An efficient shock-capturing central-type scheme for multidimensional
  relativistic flows {II}. magnetohydrodynamics.
\newblock {\em Astronomy and Astrophysics}, 400:397--413, 2003.

\bibitem{Zanotti2011}
O.~{Zanotti}, C.~{Roedig}, L.~{Rezzolla}, and L.~{Del Zanna}.
\newblock {General relativistic radiation hydrodynamics of accretion flows - I.
  Bondi-Hoyle accretion}.
\newblock {\em Mon. Not. Roy. Astr. Soc.}, 417:2899--2915, November 2011.

\end{thebibliography}
\bibliographystyle{plain}

\end{document}